\documentclass[11pt,oneside]{amsart}

\usepackage{tikz-cd}

\usepackage{newtxtext}
\usepackage{newtxmath}
\usepackage{fix-cm}

\DeclareMathAlphabet{\mathcal}{OMS}{cmsy}{m}{n}

\usepackage[arrow,curve,matrix,tips,2cell]{xy}

\usepackage[hmargin=2cm,vmargin=1.5cm]{geometry}

\usepackage{bm}

\usepackage{fancyhdr}

\usepackage[obeyspaces]{url}

\usepackage{amsmath}
\usepackage{amscd}
\usepackage{amssymb}
\usepackage{latexsym}
\usepackage{url}

\usepackage{graphicx}

\usepackage{adjustbox}

\newtheorem{theorem}{Theorem}[section]
\newtheorem*{theorem*}{Theorem}
\newtheorem{lemma}[theorem]{Lemma}
\newtheorem*{lemma*}{Lemma}
\newtheorem{corollary}[theorem]{Corollary}
\newtheorem{proposition}[theorem]{Proposition}

\newtheorem{remark}[theorem]{Remark}
\newtheorem{definition}[theorem]{Definition}
\newtheorem*{definition*}{Definition}

\newtheorem{question}[theorem]{Question}
\newtheorem*{question*}{Question}

\newtheorem{example}[theorem]{Example}
\newtheorem{examples}[theorem]{Examples}

\makeatletter
\def\revddots{\mathinner{\mkern1mu\raise\p@
\vbox{\kern7\p@\hbox{.}}\mkern2mu
\raise4\p@\hbox{.}\mkern2mu\raise7\p@\hbox{.}\mkern1mu}}
\makeatother 
\newcommand{\bgl}{\begin{equation}} 
\newcommand{\egl}{\end{equation}}
\newcommand{\bgloz}{\begin{equation*}} 
\newcommand{\egloz}{\end{equation*}}
\newcommand{\bgln}{\begin{eqnarray}} 
\newcommand{\egln}{\end{eqnarray}}
\newcommand{\bglnoz}{\begin{eqnarray*}} 
\newcommand{\eglnoz}{\end{eqnarray*}}
\newcommand{\btheo}{\begin{theorem}}
\newcommand{\etheo}{\end{theorem}}
\newcommand{\btheooz}{\begin{theorem*}}
\newcommand{\etheooz}{\end{theorem*}}

\newcommand{\blemma}{\begin{lemma}}
\newcommand{\elemma}{\end{lemma}}
\newcommand{\blemmaoz}{\begin{lemma*}}
\newcommand{\elemmaoz}{\end{lemma*}}
\newcommand{\bproof}{\begin{proof}}
\newcommand{\eproof}{\end{proof}}
\newcommand{\bbew}{\begin{beweis}}
\newcommand{\ebew}{\end{beweis}}
\newcommand{\bremark}{\begin{remark}\em}
\newcommand{\eremark}{\end{remark}}
\newcommand{\bdefin}{\begin{definition}}
\newcommand{\edefin}{\end{definition}}
\newcommand{\bdefinoz}{\begin{definition*}}
\newcommand{\edefinoz}{\end{definition*}}
\newcommand{\bex}{\begin{example}}
\newcommand{\eex}{\end{example}}
\newcommand{\bexs}{\begin{examples}}
\newcommand{\eexs}{\end{examples}}

\newcommand{\bprop}{\begin{proposition}}
\newcommand{\eprop}{\end{proposition}}
\newcommand{\bcor}{\begin{corollary}}
\newcommand{\ecor}{\end{corollary}}
\newcommand{\bfa}{\begin{cases}} 
\newcommand{\efa}{\end{cases}}

\newcommand{\bquestion}{\begin{question}}
\newcommand{\equestion}{\end{question}}
\newcommand{\bquestionoz}{\begin{question*}}
\newcommand{\equestionoz}{\end{question*}}
\newcommand{\cA}{\mathcal A}

\newcommand{\cD}{\mathcal D}
\newcommand{\cE}{\mathcal E}

\newcommand{\cI}{\mathcal I}
\newcommand{\cJ}{\mathcal J}
\newcommand{\cK}{\mathcal K}
\newcommand{\cL}{\mathcal L}
\newcommand{\cM}{\mathcal M}

\newcommand{\cO}{\mathcal O}
\newcommand{\cP}{\mathcal P}

\newcommand{\cS}{\mathcal S}
\newcommand{\cT}{\mathcal T}
\newcommand{\cU}{\mathcal U}
\newcommand{\cV}{\mathcal V}

\def\Cz{\mathbb{C}}

\def\Rz{\mathbb{R}}

\def\Zz{\mathbb{Z}}

\def\1z{\mathbb{1}}
\newcommand{\fA}{\mathfrak A}

\newcommand{\fD}{\mathfrak D}
\newcommand{\fE}{\mathfrak E}

\newcommand{\fR}{\mathfrak R}

\newcommand{\mfa}{\mathfrak a}
\newcommand{\fk}{\mathfrak k}
\newcommand{\fm}{\mathfrak m}

\newcommand{\mfp}{\mathfrak p}

\newcommand{\fr}{\mathfrak r}
\newcommand{\an}[1]{``#1''} 
\newcommand{\ti}{\tilde}

\newcommand{\lori}{\longrightarrow}

\newcommand{\ma}{\mapsto} 
\newcommand{\onto}{\twoheadrightarrow} 
\newcommand{\into}{\hookrightarrow} 

\newcommand{\ve}{\varepsilon}

\def\SEMI{\mbox{$\times\kern-2pt\vrule height5pt width.6pt \kern3pt $}}

\newcommand{\Spec}{{\rm Spec\,}} 

\newcommand{\id}{{\rm id}}
\newcommand{\alg}{{\rm alg}}

\newcommand{\Ad}{{\rm Ad\,}}

\newcommand{\reg}{^\times} 
\newcommand{\lspan}{{\rm span}} 
\newcommand{\clspan}{\overline{\lspan}} 
\newcommand{\abs}[1]{\lvert#1\rvert} 
\newcommand{\norm}[1]{\left\|#1\right\|} 
\newcommand{\defeq}{\mathrel{:=}} 

\newcommand{\dop}{\text{: }} 
\newcommand{\ilim}{\varinjlim} 
\newcommand{\supp}{{\rm supp}} 
\newcommand{\res}{{\rm res}}
\newcommand{\lge}{\left\{} 
\newcommand{\rge}{\right\}} 
\newcommand{\lru}{\left(} 
\newcommand{\rru}{\right)} 
\newcommand{\lsp}{\left\langle} 
\newcommand{\rsp}{\right\rangle} 
\newcommand{\rukl}[1]{\lru #1 \rru} 
\newcommand{\gekl}[1]{\lge #1 \rge} 
\newcommand{\spkl}[1]{\lsp #1 \rsp} 
\newcommand{\menge}[2]{\gekl{ #1 \dop #2 }} 
\def\bf1{\mathbf{1}}
\begin{document}

\title{K-theory for semigroup C*-algebras and partial crossed products}

\thispagestyle{fancy}

\author{Xin Li}

\address{Xin Li, School of Mathematics and Statistics, University of Glasgow, University Place, Glasgow G12 8QQ, United Kingdom}
\email{Xin.Li@glasgow.ac.uk}

\subjclass[2010]{Primary 46L80, 46L05; Secondary 20M18, 20Mxx}

\thanks{This project has received funding from the European Research Council (ERC) under the European Union's Horizon 2020 research
and innovation programme (grant agreement No. 817597).}

\begin{abstract}
Using the Baum-Connes conjecture with coefficients, we develop a K-theory formula for reduced C*-algebras of strongly $0$-$E$-unitary inverse semigroups, or equivalently, for certain reduced partial crossed products. In the case of semigroup C*-algebras, we obtain a generalization of previous K-theory results of Cuntz, Echterhoff and the author without having to assume the Toeplitz condition. As applications, we discuss semigroup C*-algebras of Artin monoids, Baumslag-Solitar monoids, one-relator monoids, C*-algebras generated by right regular representations of semigroups from number theory, and C*-algebras of inverse semigroups arising in the context of tilings.
\end{abstract}

\maketitle

\setlength{\parindent}{0cm} \setlength{\parskip}{0.5cm}

\section{Introduction}

Many prominent classes of C*-algebras, which played an important role in the development of the subject, are generated by partial isometries, for example AF algebras \cite{Bra,Ell}, Cuntz-Krieger algebras \cite{CK}, graph algebras \cite{Rae}, tiling C*-algebras \cite{Kel95} or semigroup C*-algebras \cite{CELY}, to mention just a few. The notions of inverse semigroups and inverse semigroup C*-algebras \cite{Pat,Ex08} provide a natural and powerful framework to study these C*-algebras and their properties. Another very general and powerful concept is given by partial dynamical systems and the corresponding partial crossed product construction (see \cite{Ex17}), which highlights the underlying dynamics at the heart of many C*-algebra constructions, including the ones mentioned above. A connection between inverse semigroups and partial dynamical systems, in particular in view of the corresponding C*-algebra constructions, is provided by the notion of strongly $0$-$E$-unitary inverse semigroups, also called strongly $E^*$-unitary inverse semigroups (see \cite{Law,BFFG,MS}).

The goal of this paper is to compute K-theory for reduced C*-algebras of strongly $0$-E-unitary inverse semigroups, i.e., C*-algebras generated by left regular representations of these inverse semigroups. This is a very natural task, given the important role K-theory plays in classification of C*-algebras \cite{Ror} and in extracting interesting information from C*-algebra constructions. In many cases, this task is challenging, but at the same time also rewarding as computing K-theory often reveals much of the inner structure of the C*-algebras in question.

The main result of this paper can also be reformulated in terms of partial dynamical systems. In this context, it provides a K-theory formula for reduced crossed products of partial dynamical systems which admit an invariant regular basis of compact open sets. This is the analogue of the corresponding notion for global dynamical systems, as introduced in \cite[\S~6]{CEL1} and \cite[\S~2]{CEL2}. Thus the K-theory formula derived in this paper generalizes the results in \cite{CEL1,CEL2} from global dynamical systems to partial dynamical systems. From the perspective of inverse semigroups, the main theorem of this paper generalizes the K-theory result in \cite{Nor15}, which applied to inverse semigroups which are not only strongly $0$-$E$-unitary but satisfy the stronger property that they are $0$-$F$-inverse semigroups and admit a partial homomorphism to a group which is injective on maximal elements. For example, while it is not clear whether the main K-theoretic result in \cite{Nor15} applies to all inverse semigroups attached to tilings, we can now compute K-theory for the reduced C*-algebras of all such tiling inverse semigroups.

The main motivation for the work in \cite{CEL1,CEL2} was to compute K-theory for semigroup C*-algebras, i.e., C*-algebras generated by left regular representations of left-cancellative semigroups. This class of C*-algebras also provides the main motivation for this work. The K-theory formula developed in this paper applies to all semigroups which satisfy the independence condition from \cite{Li12} and embed into a group which satisfies the Baum-Connes conjecture with coefficients \cite{BCH,Val}. In other words, the reversibiliby condition (which is part of the Ore condition) in \cite{CEL1} and -- even better -- the Toeplitz condition in \cite{CEL2} are no longer needed. As mentioned in \cite[\S~5.11]{CELY}, such a more general formula would be interesting as it would open the way for new applications. For example, we can now compute K-theory for the semigroup C*-algebras of all Baumslag-Solitar monoids, while the Baumslag-Solitar monoids given by presentations of the form $\spkl{a,b \ \vert \ a = b^l a b^k}^+$ with $k > 1$ were not accessible previously as they do not satisfy the Toeplitz condition, no matter which group embedding we choose. Other examples where the Toeplitz condition has not been verified but where we can now nevertheless compute K-theory include semigroup C*-algebras of Artin monoids (as long as the corresponding Artin groups satisfy the Baum-Connes conjecture with coefficients) and certain one-relator monoids as well as C*-algebras generated by right regular representations of $ax+b$-type semigroups of number-theoretic origin.

Let us now present the main results of this paper. 
\btheo
\label{thm:MainInvSgp}
Let $S$ be a countable inverse semigroup with semilattice of idempotents $E$. Assume that $S$ admits an idempotent pure partial homomorphism to a group which satisfies the Baum-Connes conjecture with coefficients. Then the K-theory of the reduced C*-algebra of $S$ is given by
$$
  K_*(C^*_{\lambda}(S)) \cong \bigoplus_{[d] \in S \backslash E\reg} K_*(C^*_{\lambda}(S_d)).
$$
\etheo
\setlength{\parindent}{0cm} \setlength{\parskip}{0cm}

Here $S \backslash E\reg$ denotes the set of orbits under the natural action of $S$ on the non-zero elements $E\reg$ of $E$, and $S_d = \menge{s \in S}{s^{-1}s = ss^{-1} = d}$. Note that $S_d$ is a group with identity $d$. It is worth mentioning that we are free to choose an idempotent pure partial homomorphism, in particular the target group; it does not have to be the universal one. This degree of flexibility turns out to be very useful in order to check the hypothesis of Theorem~\ref{thm:MainInvSgp}, as we will see in concrete examples.
\setlength{\parindent}{0cm} \setlength{\parskip}{0.5cm}

An alternative reformulation of Theorem~\ref{thm:MainInvSgp} in terms of partial dynamical systems reads as follows: Let $G$ be a countable discrete group, $X$ a second countable totally disconnected locally compact Hausdorff space, and $G \curvearrowright X$ a partial dynamical system, given by $U_{g^{-1}} \to U_g, \, x \ma g.x$. A family $\cV$ of compact open subsets of $X$ is called a $G$-invariant regular basis for the compact open subsets of $X$ if for all $g \in G$, $\cV_{g^{-1}} \defeq \menge{V \in \cV}{V \subseteq U_{g^{-1}}}$ is a regular basis for the compact open subsets of $U_{g^{-1}}$, and we have $g.\cV_{g^{-1}} = \cV_g$ for all $g \in G$. 

\btheo
\label{thm:MainPartialCroPro}
Assume that $G \curvearrowright X$ admits a $G$-invariant regular basis $\cV$ for the compact open subsets of $X$ and that $G$ satisfies the Baum-Connes conjecture with coefficients. Then the K-theory of the reduced partial crossed product of $G \curvearrowright X$ is given by
$$
  K_*(C_0(X) \rtimes_r G) \cong \bigoplus_{[V] \in G \backslash \cV\reg} K_*(C^*_{\lambda}(G_V)).
$$
\etheo
\setlength{\parindent}{0cm} \setlength{\parskip}{0cm}

Here $G \backslash \cV\reg$ denotes the set of orbits under the $G$-action on the non-empty elements $\cV\reg$ of $\cV$, and $G_V = \menge{g \in G}{g.V = V}$.
\setlength{\parindent}{0cm} \setlength{\parskip}{0.5cm}

Applied to the case of semigroup C*-algebras, we obtain the following corollary of Theorem~\ref{thm:MainInvSgp}:
\bcor
\label{cor:MainP}
Let $P$ be a countable subsemigroup of a group which satisfies the Baum-Connes conjecture with coefficients. Assume that $P$ satisfies the independence condition. Then the K-theory of the semigroup C*-algebra of $P$ is given by
$$
  K_*(C^*_{\lambda}(P)) \cong \bigoplus_{[X] \in P \backslash \cJ_P\reg} K_*(C^*_{\lambda}(P_X)).
$$
\ecor
\setlength{\parindent}{0cm} \setlength{\parskip}{0cm}

Here $\cJ_P\reg$ is the set of non-empty constructible right ideals of $P$, as introduced in \cite{Li12}, $P \backslash \cJ_P\reg$ is basically the set of orbits for the natural $P$-action on $\cJ_P\reg$ and more precisely the set of equivalence classes of the equivalence relation on $\cJ_P\reg$ generated by $X \sim pX = \menge{px}{x \in X}$ for all $X \in \cJ_P\reg$ and $p \in P$, and $P_X$ is the group of bijections $X \to X$ which can be expressed as compositions of finitely many maps, each of which given by left multiplication by a fixed semigroup element or the (set-theoretical) inverse of such a left multiplication map. 
\setlength{\parindent}{0cm} \setlength{\parskip}{0.5cm}

Note that for every countable subsemigroup of a group which satisfies the Baum-Connes conjecture with coefficients, Theorem~\ref{thm:MainInvSgp} always gives us -- no matter whether $P$ satisfies the independence condition or not -- a K-theory formula for the reduced C*-algebra of the left inverse hull $I_l(P)$ of $P$, i.e., the inverse semigroup of partial bijections of $P$ generated by the left multiplication maps $P \to pP, \, x \ma px$, for $p \in P$.

In the particular case of right LCM monoids (monoids for which all constructible right ideals are principal), the K-theory formula in Corollary~\ref{cor:P} specializes as follows:
\bcor
\label{cor:MainLCM}
Let $P$ be a countable right LCM monoid, and assume that $P$ embeds into a group which satisfies the Baum-Connes conjecture with coefficients. Then the K-theory of the semigroup C*-algebra of $P$ is given by
$$
  K_*(C^*_{\lambda}(P)) \cong K_*(C^*_{\lambda}(P^*)).
$$
\ecor
\setlength{\parindent}{0cm} \setlength{\parskip}{0cm}

Here $P^*$ is the group of invertible elements in $P$.
\setlength{\parindent}{0cm} \setlength{\parskip}{0.5cm}

In all the above-mentioned results, the K-theory isomorphisms are implemented by concrete homomorphisms. If the groups in question satisfy the strong Baum-Connes conjecture in the sense of \cite[Definition~3.4.17]{CELY}, then we can actually replace K-theory isomorphism by KK-equivalence. Moreover, it is worth pointing out that the Baum-Connes conjecture is not really needed for all coefficients, but only for particular coefficients, as we will see. In addition, it is possible to add a coefficient algebra together with an automorphic action as in \cite{CEL2}.

As for the proofs of the main results, the main innovation is to utilize the Morita enveloping action \cite{Aba,Ex17}, which allows us to identify up to Morita equivalence -- and hence in K-theory -- a partial crossed product with an associated global crossed product. However, even if we start with a partial dynamical system on a locally compact Hausdorff space, i.e., on a commutative C*-algebra, the Morita enveloping action will not be on a commutative C*-algebra anymore. This means that the K-theory results in \cite{CEL1,CEL2} do not apply. Nevertheless, it turns out that using the Baum-Connes conjecture with coefficients in the form of the Going-Down principle \cite{CEO, ENO, MN, ELPW}, it is possible to compute K-theory for the crossed product of the Morita enveloping action by comparing it to a discrete version. This part of the proof is in spirit very similar to the strategy in \cite{CEL1,CEL2}. However, even in the setting of \cite{CEL1,CEL2}, i.e., for global dynamical systems admitting an invariant regular basis of compact open sets, the proof in the present paper differs from the one in \cite{CEL1,CEL2} because the Morita enveloping action is always (except in trivial cases) defined on a noncommutative C*-algebra. For global dynamical systems, the route in \cite{CEL1,CEL2} is more direct as it is not necessary to pass over to the Morita enveloping action. But for partial dynamical systems which do not admit a globalization in the topological setting, i.e., for which the underlying space of the (topological) enveloping action is not Hausdorff (see \cite[\S~1]{Aba} for examples), the passage to the Morita enveloping action is absolutely crucial as it allows us to apply the Going-Down principle.

The paper is structured as follows: We start with a preliminary section (\S~\ref{s:Pre}) on strongly $0$-$E$-unitary inverse semigroups, partial dynamical systems, Morita enveloping actions and the Going-Down principle. In the main section, \S~\ref{s:KMoritaEnv}, we introduce the discrete version of the crossed product of the Morita enveloping action, construct an KK-element (both in \S~\ref{ss:DiscreteV}) and show that this element implements an KK-equivalence between the crossed product of the Morita enveloping action and its discrete version. This uses the Going-Down principle as well as inductive limit decompositions (\S~\ref{ss:IndLimDec}) and a careful analysis of the finite-dimensional case, where our KK-element can be described by a finite-dimensional matrix, and showing invertibility boils down to decomposing this matrix into the sum of the identity matrix and a nilpotent matrix (\S~\ref{ss:KKEqFiniteSubg}). In \S~\ref{ss:HomIndKKEq}, we construct explicit homomorphisms which induce the K-theory isomorphisms in Theorems~\ref{thm:MainInvSgp}, \ref{thm:MainPartialCroPro} and Corollaries~\ref{cor:MainP}, \ref{cor:MainLCM}. We complete the proof of our main theorems in \S~\ref{ss:Pfs}. Finally, we present applications in \S~\ref{s:App}. We compute K-theory for semigroup C*-algebras of Artin monoids, Baumslag-Solitar monoids and one-relator monoids. In the later case, our K-theory formula leads to a classification result for semigroup C*-algebras. Moreover, we compute K-theory for the C*-algebras generated by right regular representations of $ax+b$-type semigroups of number-theoretic origin. As a consequence, we deduce that for this class of semigroups, the C*-algebras generated by their left and right regular representations are KK-equivalent. This is an interesting phenomenon which already appeared in \cite[\S~6]{CEL2} and \cite[\S~4]{Li16} (see also the discussion in \cite[\S~5.11]{CELY}). Finally, we compute K-theory for reduced C*-algebras of tiling and point-set inverse semigroups as well as of another class of closely related inverse semigroups.

\section{Preliminaries}
\label{s:Pre}

Let us first recall some basic structures and constructions which will play an important role in this paper.

\subsection{Inverse semigroups}
\label{ss:InvSgp}

The interested reader may consult \cite{Pat,Ex08,Ex17,CELY} for more details and references concerning the contents of this subsection.

\bdefin
An inverse semigroup is a semigroup $S$ with the property that for every $s \in S$, there exists a unique element in $S$, denoted by $s^{-1}$, such that $s = s s^{-1} s$ and $s^{-1} = s^{-1} s s^{-1}$.
\setlength{\parindent}{0.5cm} \setlength{\parskip}{0cm}

An inverse semigroup with zero is an inverse semigroup $S$ together with a distinguished element $0 \in S$ such that $0 \cdot s = 0 = s \cdot 0$ for all $s \in S$.
\edefin
\setlength{\parindent}{0cm} \setlength{\parskip}{0cm}

One way to think about inverse semigroups is to view them as semigroups of partial bijections on a given set. The multiplication is then given by composition of partial bijections (where the domains and ranges have to be adjusted accordingly).
\setlength{\parindent}{0cm} \setlength{\parskip}{0.5cm}

Note that, given an inverse semigroup $S$, we can always construct an inverse semigroup with zero by adding $0$: As the underlying set, consider $S \cup \gekl{0}$, keep the multiplication on $S$, and define $0 \cdot s \defeq 0$ and $s \cdot 0 \defeq 0$ for all $s \in S$ as well as $0 \cdot 0 \defeq 0$. Therefore, in this paper, we will only consider inverse semigroups with zero. For the sake of brevity, we will simply write \an{inverse semigroup} instead of \an{inverse semigroup with zero}.

For the study of semigroup C*-algebras, the following examples of inverse semigroups are important: Let $P$ be a left cancellative semigroup. Let $I_l(P)$ be the smallest inverse semigroup of partial bijections on $P$ containing all partial bijections of the form $P \to pP, \, x \ma px$ which are given by left multiplication by a semigroup element $p \in P$ as well as the partial bijection $\emptyset \to \emptyset$. The later element is denoted by $0$ and is the zero element in $I_l(P)$. We call $I_l(P)$ the left inverse hull of $P$.

The following is an important piece of structure in inverse semigroups:
\bdefin
The semilattice of idempotents of an inverse semigroup $S$ is given by $E \defeq \menge{s^{-1}s}{s \in S} = \menge{s s^{-1}}{s \in S} = \menge{e \in S}{e = e^2}$.
\setlength{\parindent}{0.5cm} \setlength{\parskip}{0cm}

Given $d, e \in E$, we write $d \leq e$ if $d = d \cdot e$.
\edefin
\setlength{\parindent}{0cm} \setlength{\parskip}{0cm}

If we view inverse semigroups as sets of partial bijections, then the semilattice of idempotents can be identified with the domains and ranges of the partial bijections. For example, if our inverse semigroup is given by the left inverse hull $I_l(P)$ as above, then its semilattice of idempotents can be identified with the semilattice of subsets of $P$ which can be obtained from $P$ be finitely many operations, each of which is given by left multiplication by a fixed semigroup element of $P$ or by taking the pre-image under left multiplication by a fixed semigroup element. These subsets are called constructible right ideals of $P$, and we write $\cJ_P \defeq E(I_l(P))$.
\setlength{\parindent}{0cm} \setlength{\parskip}{0.5cm}

Let us now recall the notion of strongly $0$-$E$-unitary inverse semigroups from \cite{Law,BFFG}. Given an inverse semigroup $S$ with zero $0$, set $S\reg \defeq S \setminus \gekl{0}$ and $E\reg \defeq E \setminus \gekl{0}$. Let $G$ be a group with identity $1$.
\bdefin
A map $\sigma: \: S\reg \to G$ is called a partial homomorphism if $\sigma(st) = \sigma(s) \sigma(t)$ for all $s,t \in S\reg$ with $st \in S\reg$. 
\setlength{\parindent}{0.5cm} \setlength{\parskip}{0cm}

A map $\sigma: \: S\reg \to G$ is called idempotent pure if $\sigma^{-1}(1) = E\reg$.

We call $S$ strongly $0$-$E$-unitary if it admits an idempotent pure partial homomorphism to a group.
\edefin
\setlength{\parindent}{0cm} \setlength{\parskip}{0.5cm}

We recall the following useful observation.
\blemma[Lemma~5.5.7 in \cite{CELY}]
\label{lem:s=ge}
Let $S$ be an inverse semigroup with an idempotent pure partial homomorphism $\sigma: \: S\reg \to G$. Let $s,t$ be elements of $S\reg$. If $s^{-1} s = t^{-1} t$ and $\sigma(s) = \sigma(t)$, then $s = t$.
\elemma

For example, suppose $P$ is a subsemigroup of a group $G$. Then for every partial bijection $s \in I_l(P)\reg$, there is a unique element $\sigma(s) \in G$ such that $s(x) = \sigma(s) \cdot x$ for all $x$ in the domain of $s$. This allows us to define a map $\sigma: \: I_l(P)\reg \to G$ by sending $s \in I_l(P)\reg$ to $\sigma(s)$, and it is easy to see that this is an idempotent pure partial homomorphism.

Let us now construct reduced C*-algebras of inverse semigroups. We start with left regular representations. Let $S$ be an inverse semigroup. For $s \in S$, define $\lambda_s: \: \ell^2 S\reg \to \ell^2 S\reg$ by $\lambda_s(\delta_x) \defeq \delta_{sx}$ if $s^{-1}s \geq xx^{-1}$ and $\lambda_s(\delta_x) \defeq 0$ otherwise. Here $\delta_x$ is the element of $\ell^2 S\reg$ given by $\delta_x(y) = 1$ if $x=y$ and $\delta_x(y) = 0$ if $x \neq y$. Now we consider the C*-algebra generated by the left regular representation $\lambda$.
\bdefin
The reduced C*-algebra of $S$ is given by $C^*_{\lambda}(S) \defeq C^*(\menge{\lambda_s}{s \in S}) \subseteq \cL(\ell^2 S\reg)$.
\edefin
\setlength{\parindent}{0cm} \setlength{\parskip}{0cm}

Note that $C^*_{\lambda}(S)$ contains $C^*(E) \defeq C^*(\menge{\lambda_e}{e \in E}$ as a commutative sub-C*-algebra. 
\setlength{\parindent}{0cm} \setlength{\parskip}{0.5cm}

If we start with an inverse semigroup and adjoin a zero as we did above, then it is easy to see that the reduced C*-algebra of the original inverse semigroup (as defined in \cite{Pat}) coincides with our reduced C*-algebra of the new inverse semigroup with zero. 

Moreover, it is worth pointing out that there is a full version as well, i.e., the C*-algebra which is universal for all representations of our inverse semigroup as partial isometries. We will not need the full version in this paper. For semilattices, however, the reduced and full versions always coincide, which is why we simply write $C^*(E)$. For simplicity, given $e \in E$, we write $e$ for the element $\lambda_e \in C^*(E)$. Furthermore, there are tight versions of the reduced and full C*-algebras attached to inverse semigroups, which are given by natural quotients.

Let us now recall the construction of reduced semigroup C*-algebras, i.e., C*-algebras generated by left regular representations of left-cancellative semigroups. Given such a semigroup $P$, define for every $p \in P$ the isometry $V_p: \: \ell^2P \to \ell^2P$ by $V_p(\delta_x) \defeq \delta_{px}$. 
\bdefin
The reduced semigroup C*-algebra of $P$ is given by $C^*_{\lambda}(P) \defeq C^*(\menge{V_p}{p \in P}) \subseteq \cL(\ell^2P)$.
\edefin

Let us now compare $C^*_{\lambda}(P)$ with $C^*_{\lambda}(I_l(P))$. We denote by $p$ the partial bijection $P \cong pP, \, x \ma px$, which lies in $I_l(P)$. As explained in \cite[\S~3.2]{Nor14} or \cite[Lemma~5.6.11]{CELY}, there always exists a surjective homomorphism 
\begin{equation}
\label{e:ItoP}
  C^*_{\lambda}(I_l(P)) \onto C^*_{\lambda}(P) \ {\rm sending} \ \lambda_p \ {\rm to} \ V_p.
\end{equation}
To explain when this homomorphism is an isomorphism, we need the independence condition from \cite{Li12}.

\bdefin
\label{def:Ind}
A left cancellative semigroup $P$ is said to satisfy the independence condition if for every $X \in \cJ_P$ and all $X_1, \dotsc, X_n \in \cJ_P$, $X = \bigcup_{i=1}^n X_i$ implies that $X = X_i$ for some $1 \leq i \leq n$.
\edefin

The following result has been observed in \cite{Nor14} (see also \cite[Proposition~5.6.37]{CELY}):
\bprop
\label{prop:P=S}
Let $P$ be a subsemigroup of a group. Then the map in \eqref{e:ItoP} is an isomorphism if and only if $P$ satisfies the independence condition.
\eprop

\subsection{Partial dynamical systems}
\label{ss:PDS}

Let us now recall the notion of partial dynamical systems and explain the connection to inverse semigroups. We refer to \cite{McCl,Aba,Ex08,Ex17,CELY} for more details and references concerning the contents of this subsection.

All the groups in this paper are discrete (but see the beginning of \S~\ref{s:KMoritaEnv} for a discussion about this). Let $G$ be a group with identity $1$ and $X$ a topological space (which for us will always be locally compact and Hausdorff).
\bdefin
A partial dynamical system $G \curvearrowright X$ is given by homeomorphisms $\alpha_g: \: U_{g^{-1}} \to U_g, \, x \ma g.x$ for all $g \in G$, where $U_g$ are open subsets of $X$, such that $\alpha_1 = \id_X$ (in particular, $U_1 = X$) and for all $g, h \in G$, we have $h.(U_{(gh)^{-1}} \cap U_{h^{-1}}) = U_h \cap U_{g^{-1}}$ and $(gh).x = g.(h.x)$ for all $x \in U_{(gh)^{-1}} \cap U_{h^{-1}}$.
\edefin
The dual partial action of $\alpha$ is given by $\alpha_g^*: \: C_0(U_{g^{-1}}) \to C_0(U_g), \, f \ma f(g^{-1}.\sqcup)$. $\alpha^*$ is a partial action of $G$ on $C_0(X)$ in the sense of \cite{McCl}. Now recall that the reduced crossed product $C_0(X) \rtimes_{\alpha^*,r} G$ (or just $C_0(X) \rtimes_r G$ if $\alpha^*$ is understood) is given by the completion of $C_0(X) \rtimes^{\alg} G \defeq \menge{\sum_g f_g \delta_g \in C_c(G,C_0(X))}{f_g \in C_0(U_g)}$ or $C_0(X) \rtimes^{\ell^1} G \defeq \menge{\sum_g f_g \delta_g \in \ell^1(G,C_0(X))}{f_g \in C_0(U_g)}$ under a norm induced by a concrete representation generalizing the construction of reduced crossed products for global dynamical systems. Here the *-algebra structure on $C_0(X) \rtimes^{\alg} G \defeq \menge{\sum_g f_g \delta_g \in C_c(G,C_0(X))}{f_g \in C_0(U_g)}$ or $C_0(X) \rtimes^{\ell^1} G \defeq \menge{\sum_g f_g \delta_g \in \ell^1(G,C_0(X))}{f_g \in C_0(U_g)}$ is given by $\rukl{\sum_g f_g \delta_g} \cdot \rukl{\sum_h f'_h \delta_h} \defeq \sum_{g,h} \alpha^*_g(\alpha^*_{g^{-1}} (f_g) f'_h) \delta_{gh}$ as multiplication and $\rukl{\sum_g f_g \delta_g}^* \defeq \sum_g \alpha^*_g(f^*_{g^{-1}}) \delta_g$ as involution. We refer for more details to \cite{McCl} or \cite[\S~5.5.2]{CELY}. Again, there is also a notion of full crossed products for partial dynamical systems, which we will not need.

\bex
\label{ex:S-pDS}
Suppose that $S$ is an inverse semigroup with an idempotent pure partial homomorphism $\sigma: \: S\reg \to G$ to a group $G$. Let $E$ be the semilattice of idempotents of $S$. Define $\widehat{E}$ to be the space of non-zero homomorphisms $E \to \gekl{0,1}$ sending $0 \in E$ to $0$. We equip $\widehat{E}$ with the topology of pointwise convergence. It is easy to see that $\widehat{E}$ is canonically homeomorphic to $\Spec(C^*(E))$. Let us now construct a partial dynamical system $G \curvearrowright \widehat{E}$ as follows: The dual partial action $\alpha^*$ of $G$ on $D \defeq C^*(E) \cong C_0(\widehat{E})$ is given by $D_{g^{-1}} \defeq \clspan(\menge{s^{-1}s}{s \in S\reg, \, \sigma(s) = g})$ and $\alpha^*_g: \: D_{g^{-1}} \to D_g, \, s^{-1}s \ma ss^{-1}$. To describe the partial action $\alpha$ of $G$ on $\widehat{E}$, set $U_g \defeq \Spec(D_g) \subseteq \Spec(D) \cong \widehat{E}$. It is easy to check that $U_{g^{-1}} = \menge{\chi \in \widehat{E}}{\exists \ s \in S\reg, \, \sigma(s) = g \text{ such that } \chi(s^{-1}s) = 1}$. Given $\chi \in U_{g^{-1}}$ and $s \in S\reg$ with $\sigma(s) = g$ and $\chi(s^{-1}s) = 1$, we then have $\alpha_g(\chi)(e) = \chi(s^{-1} e s)$ for all $e \in E$.
\eex

As observed in \cite{Li17} (see also \cite[Corollary~5.5.23]{CELY}), we have an explicit isomorphism
\begin{equation}
\label{e:CS=DxG}
C^*_{\lambda}(S) \cong D \rtimes_r G, \, \lambda_s \ma (ss^{-1}) \delta_{\sigma(s)} \ {\rm for} \ s \in S\reg.
\end{equation}

Note that the partial dynamical system $G \curvearrowright D$ restricts to a partial dynamical system $G \curvearrowright E$. To see this, we need the following observation:
\blemma
\label{lem:Eg=s-1s}
$E\reg \cap D_{g^{-1}} = \menge{s^{-1}s}{s \in S\reg, \, \sigma(s) = g}$.
\elemma
\setlength{\parindent}{0cm} \setlength{\parskip}{0cm}

\bproof
\an{$\supseteq$} is clear. To prove \an{$\subseteq$}, take $e \in E\reg$ with $e \in D_{g^{-1}}$. Since $D_{g^{-1}} \defeq \clspan(\menge{s^{-1}s}{s \in S\reg, \, \sigma(s) = g})$, there must exist $d_1, \dotsc, d_n \in \menge{s^{-1}s}{s \in S\reg, \, \sigma(s) = g}$ and $\alpha_1, \dotsc, \alpha_n \in \Cz$ such that $\norm{e - \sum_i \alpha_i d_i} < \tfrac{1}{4}$. Without loss of generality we may assume that $\gekl{d_1, \dotsc, d_n}$ is multiplicatively closed. Then, by functional calculus, we deduce that $e \in C^*(\gekl{d_1, \dotsc, d_n}) = \lspan(\gekl{d_1, \dotsc, d_n})$. Hence it follows that $e = \bigvee_i e d_i$. But then \cite[Proposition~3.4]{Nor14} implies that $e = e d_i$ for some $1 \leq i \leq n$ (see also \cite[Definition~2.6 and Remark~2.7]{CEL2}). Thus $e = ed_i \in \menge{s^{-1}s}{s \in S\reg, \, \sigma(s) = g}$.
\eproof
As a consequence, it follows that $G \curvearrowright D$ restricts to the partial dynamical system $G \curvearrowright E$ given by $E_{g^{-1}} \defeq E \cap D_{g^{-1}} = \menge{s^{-1}s}{s \in S\reg, \, \sigma(s) = g}$ and $g.(s^{-1}s) = ss^{-1} \in E_g$.
\setlength{\parindent}{0cm} \setlength{\parskip}{0.5cm}

We have thus seen that given an inverse semigroup with an idempotent pure partial homomorphism to a group, we can construct a partial dynamical system of that group such that its reduced crossed product is canonically isomorphic to the reduced C*-algebra of the original inverse semigroup. But which partial dynamical systems arise in this way from inverse semigroups? This question leads to the following generalized notion of invariant regular basis of compact open subsets. The corresponding notion for global dynamical systems has been introduced in \cite{CEL1,CEL2}. 
\bdefin
\label{def:IndpDS}
Let $G \curvearrowright X$ be a partial dynamical system. A family $\cV$ of compact open subsets of $X$ is called a $G$-invariant regular basis of compact open subsets of $X$ if for all $g \in G$, $\cV_{g^{-1}} \defeq \menge{V \in \cV}{V \subseteq U_{g^{-1}}}$ is a regular basis for the compact open subsets of $U_{g^{-1}}$, in the sense of \cite[Definition~2.9]{CEL2}, and $g.\cV_{g^{-1}} = \cV_g$ (i.e., $g.V$ lies in $\cV_g$ for all $V \in \cV_{g^{-1}}$). 
\edefin
\setlength{\parindent}{0cm} \setlength{\parskip}{0cm}

Note that $\cV = \cV_1$, so that $\cV$ is a regular basis for the compact open subsets of $X$ in the sense of \cite[Definition~2.9]{CEL2}. In particular, $X$ must be locally compact Hausdorff and totally disconnected.
\setlength{\parindent}{0cm} \setlength{\parskip}{0.5cm}

Let us now explain why the partial dynamical systems appearing in Definition~\ref{def:IndpDS} are precisely those which arise from inverse semigroups as in Example~\ref{ex:S-pDS}. More precisely, the following constructions are -- up to isomorphism -- inverse to each other:
\setlength{\parindent}{0cm} \setlength{\parskip}{0cm}

\begin{itemize}
\item[($*$)] Given an inverse semigroup $S$ with an idempotent pure partial homomorphism $\sigma: \: S\reg \to G$ to a group $G$, let $G \curvearrowright X \defeq \widehat{E}$ be as in Example~\ref{ex:S-pDS}. Define $\cV \defeq \menge{\supp(e)}{e \in E}$, where we make use of the canonical isomorphism $D \defeq C^*(E) \cong C_0(X)$ under which $e$ is identified with the characteristic function $\bm{1}_{\supp(e)}$ on $\supp(e)$. 
\item[($**$)] Given a partial dynamical system $G \curvearrowright X$ and a $G$-invariant regular basis $\cV$ for the compact open subsets of $X$, construct an inverse semigroup $S$ by setting $S \defeq \menge{(g,V)}{g \in G, \, \emptyset \neq V \in \cV_{g^{-1}}} \cup \gekl{0}$ and by defining multiplication as $(h,W)(g,V) \defeq (hg,g^{-1}(W \cap g.V))$ if $W \cap g.V \neq \emptyset$ and $(h,W)(g,V) \defeq 0$ otherwise. Moreover, define $\sigma: \: S\reg \to G, \, (g,V) \ma g$.
\end{itemize}
It is easy to check that these constructions are well-defined. For instance, that $\cV_{g^{-1}}$ generates $U_{g^{-1}}$ in ($*$) follows from $D_{g^{-1}} = \clspan(E_{g^{-1}})$ and $\cV_{g^{-1}} = \menge{\supp(e)}{e \in E_{g^{-1}}}$. Moreover, it is easy to see that $S$ as defined in ($**$) is a well-defined inverse semigroup with $(g,V)^{-1} = (g^{-1},g.V)$ and semilattice of idempotents given by $\menge{(1,V)}{\emptyset \neq V \in \cV} \cup \gekl{0}$, and that $\sigma$ is an idempotent pure partial homomorphism. 
\setlength{\parindent}{0cm} \setlength{\parskip}{0.5cm}

\blemma
\label{lem:S=GX}
The constructions in ($*$) and ($**$) are inverse to each other. More precisely, the following are true:
\setlength{\parindent}{0cm} \setlength{\parskip}{0cm}

\begin{enumerate}
\item[(i)] If we start with an inverse semigroup $S$ and an idempotent pure partial homomorphism $\sigma: \: S\reg \to G$ to a group $G$, construct $G \curvearrowright \cV$ and $\cV$ as in ($*$), and then the inverse semigroup $\ti{S}$ with idempotent pure partial homomorphism $\ti{\sigma}$ as in ($**$), then we have an isomorphism $\rho: \: S \cong \ti{S}, \, S\reg \ni s \ma (\sigma(s), \supp(s^{-1}s)), \, 0 \ma 0$ such that $\ti{\sigma} \circ \rho = \sigma$.
\item[(ii)] If we start with a partial dynamical system $\alpha: \: G \curvearrowright X$ and a $G$-invariant regular basis $\cV$ for the compact open subsets of $X$, construct $S$ with semilattice of idempotents $E$ as in ($**$) and then construct a partial dynamical system $\ti{\alpha}: \: G \curvearrowright \ti{X}$ and a $G$-invariant regular basis $\ti{\cV}$ for the compact open subsets of $\ti{X}$ as in ($*$), then there is a homeomorphism $\varphi: \: \ti{X} \cong X$ sending $\chi \in \ti{X}$ to the unique point $x \in X$ with the property that for all $V \in \cV$, $x$ lies in $V$ if and only if $\chi(1,V) = 1$, and this homeomorphism $\varphi$ gives rise to a conjugacy between $\ti{\alpha}$ and $\alpha$ sending $\cV$ to $\ti{\cV}$.
\end{enumerate}
\elemma
\bproof
(i) To see that $\rho$ is a semigroup homomorphism, observe that $(\sigma(s), \supp(s^{-1}s)) (\sigma(t), \supp(t^{-1}t)) = (\sigma(s)\sigma(t), \sigma(t)^{-1}.(\supp(s^{-1}s) \cap \sigma(t).\supp(t^{-1}t)))$. So it suffices to show that $\sigma(t)^{-1}.(\supp(s^{-1}s) \cap \sigma(t).\supp(t^{-1}t)) = \supp((st)^{-1}(st))$. We have $\chi \in \supp((st)^{-1}(st))$ if and only if $\chi(t^{-1}t) = 1$ and $(\sigma(t).\chi)(s^{-1}s) = 1$ if and only if $\chi \in \sigma(t)^{-1}.(\supp(s^{-1}s) \cap \sigma(t).\supp(t^{-1}t))$. The inverse $\ti{S} \to S$ of $\rho$ is given as follows: Given $(g,V) \in \ti{S}\reg$, write $V = \supp(e)$ for some $e \in E\reg$. Then $V \subseteq U_{g^{-1}}$ implies $e \in D_{g^{-1}}$, which by Lemma~\ref{lem:Eg=s-1s} yields that there exists $s \in S\reg$ with $\sigma(s) = g$ and $e = s^{-1}s$. This element $s$ is uniquely determined by Lemma~\ref{lem:s=ge}. Now define $\ti{S}\reg \to S$ by sending $(g,V)$ to this element $s$ and $0$ to $0$. By construction, this is the inverse of $\rho$. It is clear that $\ti{\sigma} \circ \rho = \sigma$.
\setlength{\parindent}{0cm} \setlength{\parskip}{0.5cm}

(ii) By construction, we have an isomorphism $E \cong \cV$ sending $(1,V)$ to $V$ and $0$ to $\emptyset$. Thus we obtain $\ti{X} = \widehat{E} \cong \widehat{\cV} \cong X$, where the first equality is by construction (see ($*$)), the first homeomorphism is obtained by dualizing the isomorphism $E \cong \cV$, and the third homeomorphism is given by the composite $\widehat{\cV} \cong \Spec C^*(\menge{\bm{1}_V}{V \in \cV}) = \Spec C_0(X) \cong X$, where we used that $\cV$ is a regular basis for the compact open subsets of $X$. Following the definitions, it is easy to check that we indeed get the homeomorphism $\varphi$ as defined in (ii). Clearly, $\varphi(\supp(e,V)) = V$ for all $V \in \cV$, so that $\varphi$ sends $\ti{\cV}$ to $\cV$. Moreover, if $\varphi(\chi) = x$, then we have $\chi \in \ti{U}_{g^{-1}}$ if and only if there exists $V \in \cV_{g^{-1}}$ with $\chi(1,V) = 1$ if and only if there exists $V \in \cV_{g^{-1}}$ with $x \in V$ if and only if $x \in U_{g^{-1}}$, where $g \in G$ is arbitrary. This shows that $\varphi$ sends $\ti{U}_{g^{-1}}$ to $U_g$. Finally, suppose we have $\varphi(\chi) = x$ for some $\chi \in \ti{U}_{g^{-1}}$, i.e., there exists $V \in \cV_{g^{-1}}$ with $\chi(1,V) = 1$. Then $\ti{\alpha}_g(\chi)(1,W) = \chi((g,V)^{-1}(1,W)(g,V)) = \chi(1,\alpha_{g^{-1}}(W) \cap V) = 1$ if and only if $x \in \alpha_{g^{-1}}(W) \cap V$ if and only if $\alpha_g(x) \in W$. This implies $\varphi(\ti{\alpha}_g(\chi)) = \alpha_g(x)$. So this shows that $\varphi \circ \ti{\alpha}_g = \alpha_g \circ \varphi$ for all $g \in G$, as desired.
\eproof
\setlength{\parindent}{0cm} \setlength{\parskip}{0.5cm}

\bremark
Under the correspondence in Lemma~\ref{lem:S=GX}, countability of $S$ corresponds to second countability of the base space $X$ of the partial dynamical system $G \curvearrowright X$.
\eremark

\subsection{The Morita enveloping action}
\label{ss:MoritaEnv}

This subsection is based on \cite{Aba}. We follow the exposition in \cite[\S~26 -- \S~28]{Ex17}. Let $G \curvearrowright D$ be a partial dynamical system on a C*-algebra $D$ given by homeomorphisms $D_{g^{-1}} \cong D_g, \, d \ma g.d$. We start by defining $A$ as a sub-C*-algebra of $(D \rtimes_r G) \otimes \cK(\ell^2 G)$ by setting 
$$
  A \defeq \clspan(\menge{D_{\zeta^{-1} \eta} \delta_{\zeta^{-1} \eta} \otimes \ve_{\zeta,\eta}}{\zeta, \eta \in G}) \subseteq (D \rtimes_r G) \otimes \cK(\ell^2 G).
$$
Here $\ve_{\zeta,\eta}$ is the rank-one projection from $\Cz \delta_{\eta}$ to $\Cz \delta_{\zeta}$. $A$ is called the smash product in \cite[\S~26]{Ex17}. Next, we define $\fA$ as a sub-C*-algebra of $(D \rtimes_r G) \otimes \cK(\ell^2 G)$ by setting 
$$
  \fA \defeq \clspan(\menge{D_{\zeta^{-1}} \delta_{\zeta^{-1}} D_{\eta} \delta_{\eta} \otimes \ve_{\zeta,\eta}}{\zeta, \eta \in G}) \subseteq (D \rtimes_r G) \otimes \cK(\ell^2 G).
$$
$\fA$ is called the restricted smash product in \cite[\S~26]{Ex17}. By construction, $\fA$ is a sub-C*-algebra of $A$. $A$ comes with a $G$-action, where $g \in G$ acts via $\Ad(1 \otimes \lambda_g)$. Here $\lambda$ is the left regular representation of $G$ on $\ell^2G$. This $G$-action $G \curvearrowright A$ is called the Morita enveloping action of $G \curvearrowright D$. It restricts to a partial action $G \curvearrowright \fA$. The natural embedding $\fA \into A$ extends to an embedding $\fA \rtimes_r G \into A \rtimes_r G$, and it is shown in \cite[Theorem~28.8]{Ex17} that we can identify $\fA \rtimes_r G$ with a full corner in $A \rtimes_r G$ via this embedding. Hence we obtain a Morita equivalence $\fA \rtimes_r G \sim_M A \rtimes_r G$.

Furthermore, it is shown in \cite[Theorem~28.5]{Ex17} that $G \curvearrowright D$ is Morita equivalent to $G \curvearrowright \fA$, in a $G$-equivariant way. Since we will need the particular form of the imprimitivity bimodule constructed in \cite[\S~28]{Ex17}, let us recall it now. Let $M$ be the subspace $M \defeq \clspan(\menge{D_{\eta} \delta_{\eta} \otimes \ve_{1,\eta}}{\eta \in G})$ of $\fA$. Then upon identifying $D = D_1$ with $D_1 \delta_1 \otimes \ve_{1,1}$, we obtain a $D$-valued inner product by setting ${}_D \! \spkl{x,y} \defeq xy^*$. We obtain an $\fA$-valued inner product by setting $\spkl{x,y}_{\fA} \defeq x^*y$. In this way, together with the canonical module structures, $M$ becomes a $D$-$\fA$-imprimitivity bimodule. To define a partial $G$-action on $M$, define $M_g \defeq \clspan(\menge{D_g D_{g^{-1}} D_{\eta} \delta_{\eta} \otimes \ve_{1,\eta}}{\eta \in G})$ and $M_{g^{-1}} \to M_g$ by setting $g.(d \delta_h \otimes \ve_{1,\eta}) \defeq (g.d) \delta_{g \eta} \otimes \ve_{1,g \eta}$. It is shown in \cite[Theorem~28.7]{Ex17} that $M$ gives rise to a Morita equivalence $D \rtimes_r G \sim_M \fA \rtimes_r G$. Following the proof of \cite[Theorem~28.7]{Ex17}, we obtain a concrete $(D \rtimes_r G)$-$(\fA \rtimes_r G)$-imprimitivity bimodule as follows: First form the linking algebra 
$L = 
\begin{pmatrix}
D & M \\
M^* & \fA
\end{pmatrix}
$
of $M$. The partial $G$-actions on the components of $L$ give rise to a partial dynamical system $G \curvearrowright L$. Then
\begin{equation}
\label{e:MxG}
  M \rtimes_r G \defeq 
  \begin{pmatrix}
D & 0 \\
0 & 0
\end{pmatrix}
\rukl{L \rtimes_r G}
\begin{pmatrix}
0 & 0 \\
0 & \fA
\end{pmatrix}
\end{equation}
gives rise to the desired $(D \rtimes_r G)$-$(\fA \rtimes_r G)$-imprimitivity bimodule. If we denote by $\Delta_g$, $g \in G$, the canonical unitaries in (the multiplier algebra of) $\fA \rtimes_r G$, $A \rtimes_r G$ and $L \rtimes_r G$ implementing the $G$-actions, then $M \rtimes_r G$ can be alternatively described as 
$M \rtimes_r G = \clspan(\menge{x \Delta_g}{x \in M, \, g \in G})
$ (where we identify $x$ with 
$
\rukl{
\begin{smallmatrix}
0 & x \\
0 & 0
\end{smallmatrix}
}
$), with left $D \rtimes_r G$-module structure given by $(d \delta_g) \cdot (x \Delta_h) = g. ((g^{-1}.d) x) \Delta_{gh}$, $D \rtimes_r G$-valued inner product given by ${}_{D \rtimes_r G} \! \spkl{x \Delta_g, y \Delta_h} = g. ({}_D \! \spkl{g^{-1}.x,h^{-1}.y}) \delta_{gh^{-1}}$, right $\fA \rtimes_r G$-module structure given by $(x \Delta_g) \cdot (a \Delta_h) = g.( (g^{-1}.x) a) \Delta_{gh}$ and $\fA \rtimes_r G$-valued inner product given by $\spkl{x \Delta_g, y \Delta_h}_{\fA \rtimes_r G} = g^{-1}.( \spkl{x,y}_{\fA} ) \Delta_{g^{-1}h}$.

All in all, we obtain that $D \rtimes_r G \sim_M \fA \rtimes_r G \sim_M A \rtimes_r G$, so that $K_*(D \rtimes_r G) \cong K_*(A \rtimes_r G)$. Hence it suffices to compute K-theory for $A \rtimes_r G$.

\subsection{The Going-Down principle}
\label{ss:GoingDown}

Let us recall the Going-Down principle from \cite{CEO, ENO, MN, ELPW} (see also \cite[\S~3.5]{CELY}). It will play a crucial role in the proofs of our main results.

Let us focus on the version of the Going-Down principle for discrete groups. (There is also a general version for locally compact groups.) Let $G$ be a discrete group and $\cA$ and $A$ be $G$-algebras. Let $\bm{x} \in KK^G(\cA,A)$. Given a subgroup $F$ of $G$, we denote by $\res^G_F(\bm{x}) \in KK^F(\cA,A)$ the restriction of $\bm{x}$. We write $j_G(\bm{x}) \in KK(\cA \rtimes_r G,A \rtimes_r G)$ for the descent of $\bm{x}$.
\bprop
\label{prop:GD}
$ $
\setlength{\parindent}{0cm} \setlength{\parskip}{0cm}

\begin{enumerate}
\item If $G$ satisfies the Baum-Connes conjecture for $\cA$ and $A$, and for every finite subgroup $F$ of $G$, the Kasparov product with $j_F(\res^G_F(\bm{x}))$ induces an isomorphism $\sqcup \otimes j_F(\res^G_F(\bm{x})): \: K_*(\cA \rtimes_r F) \cong K_*(A \rtimes_r F)$, then the Kasparov product with $j_G(\bm{x})$ induces an isomorphism $\sqcup \otimes j_G(\bm{x}): \: K_*(\cA \rtimes_r G) \cong K_*(A \rtimes_r G)$.
\item If $G$ satisfies the strong Baum-Connes conjecture in the sense of \cite[Definition~3.4.17]{CELY}, and for every finite subgroup $F$ of $G$, $j_F(\res^G_F(\bm{x}))$ is a KK-equivalence between $\cA \rtimes_r F$ and $A \rtimes_r F$, then $j_G(\bm{x})$ is a KK-equivalence between $\cA \rtimes_r G$ and $A \rtimes_r G$.
\end{enumerate}
\eprop
\setlength{\parindent}{0cm} \setlength{\parskip}{0.5cm}

\section{K-theory for crossed products of Morita enveloping actions}
\label{s:KMoritaEnv}

Let us now prove the main results of this paper. We first explain the setting. Let us keep the convention that all our groups are discrete. Thinking about inverse semigroups, which are often viewed as discrete objects, this is a natural assumption. And as we will see, in all our examples, the groups will be discrete. However, the work in \cite{Aba} on partial dynamical systems and their Morita enveloping actions applies to all locally compact groups, and from that point of view, it is reasonable to expect that our results should extend to the setting of locally compact groups. Let $S$ be a countable inverse semigroup, $G$ a group, and $\sigma: \: S\reg \to G$ an idempotent pure partial homomorphism. By restricting to the group generated by the image of $\sigma$ if necessary, we may always assume that $G$ is countable. Let $E$ be the semilattice of idempotents in $S$. As in Example~\ref{ex:S-pDS}, we construct the partial dynamical system $G \curvearrowright D = C^*(E)$, denoted by $D_{g^{-1}} \to D_g, \, d \ma g.d$, and its restriction $G \curvearrowright E$ (introduced after Example~\ref{ex:S-pDS}), denoted by $E_{g^{-1}} \to E_g, \, e \ma g.e$. The smash product $A \subseteq (D \rtimes_r G) \otimes \cK(\ell^2 G)$ and the Morita enveloping action $G \curvearrowright A$ are defined as in \S~\ref{ss:MoritaEnv}.

\subsection{The discrete version}
\label{ss:DiscreteV}

We first define discrete versions of $A$ and $G \curvearrowright A$. 
\bdefin
Let $\cA \subseteq \cK(\ell^2(E \times G))$ be given by
$$
  \cA \defeq \clspan(\menge{\ve_{(d,\zeta), \, (\eta^{-1}\zeta.d,\eta)}}{\zeta, \eta \in G, \, d \in E_{\zeta^{-1}\eta}}).
$$
To define a $G$-action $G \curvearrowright \cA$, we let an element $g \in G$ act on $\cA$ by $\Ad(1 \otimes \lambda_g)$ under the canonical identification $\ell^2(E \times G) \cong \ell^2E \otimes \ell^2G$.
\edefin
\setlength{\parindent}{0cm} \setlength{\parskip}{0cm}

Note that both $G$-actions $G \curvearrowright A$ and $G \curvearrowright \cA$ are given by conjugation with the same unitaries. 
\setlength{\parindent}{0cm} \setlength{\parskip}{0.5cm}

To compare $G \curvearrowright A$ to $G \curvearrowright \cA$, we construct the following homomorphism.
\bdefin
We define
$$
  \Phi: \: \cA \to \cK(\ell^2E) \otimes A, \, \ve_{(d,\zeta), \, (\eta^{-1}\zeta.d,\eta)} \ma \ve_{d, \, \eta^{-1}\zeta.d} \otimes (d \delta_{\zeta^{-1}\eta} \otimes \ve_{\zeta,\eta}).
$$
\edefin
\setlength{\parindent}{0cm} \setlength{\parskip}{0cm}

It is easy to see that $\Phi$ is a well-defined homomorphism. For instance, to see that it is multiplicative, first note that $\ve_{(d,\zeta), \, (\eta^{-1}\zeta.d,\eta)} \cdot \ve_{(e,\eta'), \, (\theta^{-1}\eta'.e,\theta)} = 0$ if and only if $\eta^{-1}\zeta.d \neq e$ or $\eta \neq \eta'$ if and only if $\rukl{\ve_{d, \, \eta^{-1}\zeta.d} \otimes (d \delta_{\zeta^{-1}\eta} \otimes \ve_{\zeta,\eta})} \cdot \rukl{\ve_{e, \, \theta^{-1}\eta'.e} \otimes (e \delta_{(\eta')^{-1}\theta} \otimes \ve_{\eta',\theta})} = 0$. If $\eta^{-1}\zeta.d = e$ and $\eta = \eta'$, then 
$$
  \rukl{\ve_{d, \, \eta^{-1}\zeta.d} \otimes (d \delta_{\zeta^{-1}\eta} \otimes \ve_{\zeta,\eta})} \cdot \rukl{\ve_{\eta^{-1}\zeta.d, \, \theta^{-1}\zeta.d} \otimes ((\eta^{-1}\zeta.d) \delta_{\eta^{-1}\theta} \otimes \ve_{\eta,\theta})} = \ve_{d,\theta^{-1}\eta.d} \otimes (d \delta_{\zeta^{-1}\theta} \otimes \ve_{\zeta,\theta}). 
$$
This shows that $\Phi$ is multiplicative. Moreover, if we equip $\cK(\ell^2E) \otimes A$ with the $G$-action given by the tensor product of the trivial $G$-action on $\cK(\ell^2E)$ and the Morita enveloping action $G \curvearrowright A$, then it is easy to see that $\Phi$ is $G$-equivariant. Hence $\Phi$ gives rise to the element
$$
  \bm{x} \defeq [\Phi] \in KK^G(\cA, \, \cK(\ell^2E) \otimes A).
$$
\setlength{\parindent}{0cm} \setlength{\parskip}{0.5cm}

\subsection{Inductive limit decompositions}
\label{ss:IndLimDec}

Our goal is to apply the Going-Down principle to the element $\bm{x}$ we constructed. For that purpose, we have to show that for every finite subgroup $F$ of $G$, $\sqcup \otimes j_F(\res^G_F(\bm{x})): \: K_*(\cA \rtimes_r F) \to K_*((\cK(\ell^2E) \otimes A) \rtimes_r F)$ is an isomorphism. In order to reduce to finite-dimensional subalgebras, we develop inductive limit decompositions in this subsection.

Let us fix a finite subgroup $F$ of $G$. Let us also fix a finite $F$-invariant subset $\Sigma$ of $G$. We show the following three lemmas:

Suppose that we are given a family $\gekl{\fE_{\zeta,\eta}}_{\zeta,\eta \in \Sigma}$ of finite subsets $\fE_{\zeta,\eta}$ of $E_{\zeta^{-1}\eta}$. Let us construct a family $\gekl{\cE_{\zeta,\eta}}_{\zeta,\eta \in \Sigma}$ as follows: Divide $\Sigma \times \Sigma$ into pairwise disjoint subsets of the form $\menge{(\gamma \zeta, \gamma \eta)}{\gamma \in F} \cup \menge{(\gamma \eta, \gamma \zeta)}{\gamma \in F}$. Start with a pair $(\zeta,\eta)$ from one of these subsets. If $\eta^{-1}\zeta$ is of finite order, let $\cE_{\zeta,\eta}$ be the smallest sub-semilattice of $E_{\zeta^{-1}\eta}$ generated by $\spkl{\eta^{-1}\zeta}.\fE_{\zeta,\eta}$, and set $\cE_{\eta,\zeta} \defeq \cE_{\zeta,\eta}$. If $\eta^{-1}\zeta$ is not of finite order, let $\cE_{\zeta,\eta}$ be the smallest sub-semilattice of $E_{\zeta^{-1}\eta}$ generated by $\fE_{\zeta,\eta}$, and set $\cE_{\eta,\zeta} \defeq \eta^{-1} \zeta.\cE_{\zeta,\eta}$.
In either case, let $\cE_{\gamma \zeta, \gamma \eta} \defeq \cE_{\zeta,\eta}$, $\cE_{\gamma \eta, \gamma \zeta} \defeq \cE_{\eta,\zeta}$ for all $\gamma \in F$. Continue in this way for all the finite subsets whose union is $\Sigma \times \Sigma$.
\blemma
\label{lem:fE-cE}
The procedure above gives rise to a well-defined family $\gekl{\cE_{\zeta,\eta}}_{\zeta,\eta \in \Sigma}$ of finite subsemilattices $\cE_{\zeta,\eta}$ of $E_{\zeta^{-1}\eta}$ containing $\fE_{\zeta,\eta}$ which satisfies
\setlength{\parindent}{0cm} \setlength{\parskip}{0cm}

\begin{enumerate}
\item[a)] $\cE_{\gamma\zeta,\gamma\eta} = \cE_{\zeta,\eta}$ for all $\gamma \in F$ and $\zeta, \eta \in \Sigma$;
\item[b)] $\cE_{\eta,\zeta} = \eta^{-1} \zeta. \cE_{\zeta,\eta}$ for all $\zeta, \eta \in \Sigma$.
\end{enumerate}
\setlength{\parindent}{0cm} \setlength{\parskip}{0.5cm}
\elemma

To explain the next construction, we introduce the following notation: Given $\alpha, \beta \in G$, $e \in E_{\alpha^{-1}}$ and $f \in E_{\beta^{-1}}$, define $e \bullet (\alpha^{-1} f) \defeq \alpha^{-1} ((\alpha.e) f) \in E_{(\beta \alpha)^{-1}}$. Now suppose that we are given a family $\gekl{\cE_{\zeta,\eta}}_{\zeta,\eta \in \Sigma}$ of finite subsemilattices $\cE_{\zeta,\eta}$ of $E_{\zeta^{-1}\eta}$ satisfying a) and b) from Lemma~\ref{lem:fE-cE}. Let $E_{\zeta,\theta}$ be the smallest sub-semilattice of $E_{\zeta^{-1}\theta}$ generated by
$$
  \menge{
  e_0 \bullet (\zeta^{-1} \eta_1 e_1) \bullet \dotso \bullet (\zeta^{-1} \eta_{l-1} e_{l-1})
  }
  {
  \eta_0 = \zeta, \eta_1, \dotsc, \eta_{l-1} \in \Sigma, \eta_l = \theta, \, e_k \in \cE_{\eta_k,\eta_{k+1}} \ \forall \ 0 \leq k \leq l-1
  }
$$
Note that we should have written $( \dotso (e_0 \bullet (\zeta^{-1} \eta_1 e_1)) \bullet \dotso ) \bullet (\zeta^{-1} \eta_{l-1} e_{l-1})$ to be more precise, but we leave out the parentheses for the sake of readability.
\blemma
\label{lem:cE-E}
The construction above yields a family $\gekl{E_{\zeta,\eta}}_{\zeta,\eta \in \Sigma}$ of finite subsemilattices $E_{\zeta,\eta}$ of $E_{\zeta^{-1}\eta}$ containing $\cE_{\zeta,\eta}$ which satisfies a) and b) from Lemma~\ref{lem:fE-cE} as well as
\setlength{\parindent}{0cm} \setlength{\parskip}{0cm}

\begin{enumerate}
\item[c)] $E_{\zeta,\eta} \bullet (\zeta^{-1} \eta E_{\eta,\theta}) \subseteq E_{\zeta,\theta}$ for all $\zeta, \eta, \theta \in \Sigma$.
\end{enumerate}
\setlength{\parindent}{0cm} \setlength{\parskip}{0.5cm}
\elemma

\blemma
\label{lem:E-FinDimAlg}
If $\gekl{E_{\zeta,\eta}}_{\zeta,\eta \in \Sigma}$ is a family of finite subsemilattices $E_{\zeta,\eta}$ of $E_{\zeta^{-1}\eta}$ satisfying a), b) from Lemma~\ref{lem:fE-cE} and c) from Lemma~\ref{lem:cE-E}, then
\begin{equation}
\label{e:AlgE}
  \lspan(\menge{d \delta_{\zeta^{-1}\eta} \otimes \ve_{\zeta,\eta}}{d \in E_{\zeta,\eta}, \, \zeta, \eta \in \Sigma})
\end{equation}
is an $F$-invariant finite-dimensional sub-C*-algebra of $A$.
\elemma

\bproof[Proof of Lemma~\ref{lem:fE-cE}]
If $\eta^{-1}\zeta$ is of finite order, then all $\cE_{\gamma \zeta, \gamma \eta}$ coincide and are $\eta^{-1} \zeta$-invariant. This is clearly well-defined, and a) and b) are satisfied. If $\eta^{-1}\zeta$ is not of finite order, then we claim that we have $(\gamma \zeta, \gamma \eta) \neq (\gamma' \eta, \gamma' \zeta)$ for all $\gamma, \gamma' \in F$ with $\gamma \neq \gamma'$. Indeed, suppose that $\gamma \zeta = \gamma' \eta$ and $\gamma \eta = \gamma' \zeta$. Then $\zeta = \gamma^{-1} \gamma' \eta = \gamma^{-1} \gamma' \gamma^{-1} \gamma' g$ and thus $\eta^{-1} \zeta = \eta^{-1} \gamma^{-1} \gamma' \eta$, and $(\gamma^{-1} \gamma') (\gamma^{-1} \gamma') = 1$, in contradiction to our assumption that $\eta^{-1}\zeta$ is not of finite order. This shows that we can define $\cE_{\gamma \zeta, \gamma \eta} \defeq \cE_{\zeta,\eta}$, $\cE_{\gamma \eta, \gamma \zeta} \defeq \cE_{\eta,\zeta}$ for all $\gamma \in F$. Thus $\gekl{\cE_{\zeta,\eta}}_{\zeta,\eta \in \Sigma}$ is well-defined. Properties a) and b) are satisfied by construction.
\eproof

\bproof[Proof of Lemma~\ref{lem:cE-E}]
We first show that for all $\zeta, \eta \in \Sigma$, the semilattice $E_{\zeta,\eta}$ we constructed is finite. For that purpose, we set
$$
  \Pi \defeq 
  \menge{
  e_0 \bullet (\zeta^{-1} \eta_1 e_1) \bullet \dotso \bullet (\zeta^{-1} \eta_{l-1} e_{l-1})
  }
  {
  \eta_0 = \zeta, \eta_1, \dotsc, \eta_{l-1} \in \Sigma, \eta_l = \theta, \, e_k \in \cE_{\eta_k,\eta_{k+1}} \ \forall \ 0 \leq k \leq l-1
  }.
$$
In order to show that $\Pi$ is finite, we need the following observation.
\blemma
\label{lem:RedundantFactor}
If in a product 
$$
  e_0 \bullet (\zeta^{-1} \eta_1 e_1) \bullet \dotso \bullet (\zeta^{-1} \eta_{l-1} e_{l-1})
$$
two of the factors coincide, i.e., $\eta_k = \eta_{\bar{k}}$ and $e_k = e_{\bar{k}}$ for $k < \bar{k}$, then we can leave out the $\bar{k}$th factor without changing the product, i.e.,
$$
  e_0 \bullet (\zeta^{-1} \eta_1 e_1) \bullet \dotso \bullet (\zeta^{-1} \eta_{l-1} e_{l-1})
  = 
  e_0 \bullet (\zeta^{-1} \eta_1 e_1) \bullet \dotso \bullet (\zeta^{-1} \eta_{\bar{k}-1} e_{\bar{k}-1}) \bullet (\zeta^{-1} \eta_{\bar{k}+1} e_{\bar{k}+1}) \bullet \dotso \bullet (\zeta^{-1} \eta_{l-1} e_{l-1}).  
$$
\elemma
\setlength{\parindent}{0cm} \setlength{\parskip}{0cm}

\bproof[Proof of Lemma~\ref{lem:RedundantFactor}]
Let us write $\check{e} \defeq e_0 \bullet (\zeta^{-1} \eta_1 e_1) \bullet \dotso \bullet (\zeta^{-1} \eta_{k-1} e_{k-1})$. We want to show that
$$
  \check{e} \bullet (\zeta^{-1} \eta_k e_k) \bullet \dotso \bullet (\zeta^{-1} \eta_{\bar{k}-1} e_{\bar{k}-1}) \bullet (\zeta^{-1} \eta_{\bar{k}} e_{\bar{k}})
  = 
  \check{e} \bullet (\zeta^{-1} \eta_k e_k) \bullet \dotso \bullet (\zeta^{-1} \eta_{\bar{k}-1} e_{\bar{k}-1}).
$$
First, let us show that
$$
  \check{e} \bullet (\zeta^{-1} \eta_k e_k) \bullet \dotso \bullet (\zeta^{-1} \eta_{\bar{k}-1} e_{\bar{k}-1})
  \leq
  \check{e} \bullet (\zeta^{-1} \eta_k e_k)
$$
by induction on $\bar{k} - k$. The initial step $\bar{k} - k = 1$ is obvious. Now we have
$$
  \check{e} \bullet (\zeta^{-1} \eta_k e_k) \bullet \dotso \bullet (\zeta^{-1} \eta_{\bar{k}-1} e_{\bar{k}-1})
  \leq
  \check{e} \bullet (\zeta^{-1} \eta_k e_k) \bullet \dotso \bullet (\zeta^{-1} \eta_{\bar{k}-2} e_{\bar{k}-2})
  \leq
  \check{e} \bullet (\zeta^{-1} \eta_k e_k)
$$
where we used the induction hypothesis in the second inequality. Hence we have
$$
  \eta_{\bar{k}}^{-1} \zeta.(\check{e} \bullet (\zeta^{-1} \eta_k e_k) \bullet \dotso \bullet (\zeta^{-1} \eta_{\bar{k}-1} e_{\bar{k}-1}))
  \leq 
  \eta_{\bar{k}}^{-1} \zeta.(\check{e} \bullet (\zeta^{-1} \eta_k e_k))
  = (\eta_{\bar{k}}^{-1} \zeta.\check{e}) e_k
  \leq e_k = e_{\bar{k}}
$$
and thus
\begin{align*}
  &\check{e} \bullet (\zeta^{-1} \eta_k e_k) \bullet \dotso \bullet (\zeta^{-1} \eta_{\bar{k}-1} e_{\bar{k}-1}) \bullet (\zeta^{-1} \eta_{\bar{k}} e_{\bar{k}})
  =
  \zeta^{-1} \eta_{\bar{k}}.(\eta_{\bar{k}}^{-1} \zeta.(\check{e} \bullet (\zeta^{-1} \eta_k e_k) \bullet \dotso \bullet (\zeta^{-1} \eta_{\bar{k}-1} e_{\bar{k}-1})) e_{\bar{k}})\\
  =
  \ &\zeta^{-1} \eta_{\bar{k}}.(\eta_{\bar{k}}^{-1} \zeta.(\check{e} \bullet (\zeta^{-1} \eta_k e_k) \bullet \dotso \bullet (\zeta^{-1} \eta_{\bar{k}-1} e_{\bar{k}-1})))
  =
  \check{e} \bullet (\zeta^{-1} \eta_k e_k) \bullet \dotso \bullet (\zeta^{-1} \eta_{\bar{k}-1} e_{\bar{k}-1}).
\end{align*}
\eproof
Therefore, $\# \Pi \leq 2^{(\# \Sigma)^2 \cdot \# \, (\amalg_{\zeta,\eta \in \Sigma} \ \cE_{\zeta,\eta})}$, so that $\Pi$ is finite, and hence also $E_{\zeta,\eta}$.
\setlength{\parindent}{0cm} \setlength{\parskip}{0.5cm}

It remains to check a), b) and c). For a), take $\gamma \in F$. Given $\eta_0 = \zeta$, $\eta_1, \dotsc, \eta_{l-1} \in \Sigma$, $\eta_l = \theta$ and $e_k \in \cE_{\eta_k,\eta_{k+1}}$ for all $0 \leq k \leq l-1$, we have $e_k \in \cE_{\gamma \eta_k, \gamma \eta_{k+1}}$ for all $0 \leq k \leq l-1$ because we have property a) for the family $\gekl{\cE_{\zeta,\eta}}_{\zeta,\eta \in \Sigma}$, and 
$$
  e_0 \bullet (\zeta^{-1} \eta_1 e_1) \bullet \dotso \bullet (\zeta^{-1} \eta_{l-1} e_{l-1})
  =
  e_0 \bullet ((\gamma \zeta)^{-1} (\gamma \eta_1) e_1) \bullet \dotso \bullet ((\gamma \zeta)^{-1} (\gamma \eta_{l-1}) e_{l-1})
$$
which lies in $E_{\gamma \zeta, \gamma \theta}$. Hence $E_{\zeta,\theta} \subseteq E_{\gamma \zeta, \gamma \theta}$, and by symmetry, a) follows.

Next, we prove c). We need another observation.
\blemma
\label{lem:eff'}
For all $\check{e} \in E_{\theta^{-1}\zeta}$, $f \in E_{\mu^{-1} \theta}$ and $f' \in E$, we have
$$
  \check{e} \bullet ( \zeta^{-1} \theta ( f \bullet (\theta^{-1} \mu f')))
  =
  (\check{e} \bullet (\zeta^{-1} \theta f)) \bullet (\zeta^{-1} \mu f').
$$
\elemma
\setlength{\parindent}{0cm} \setlength{\parskip}{0cm}

\bproof[Proof of Lemma~\ref{lem:eff'}]
We have
\begin{align*}
  &\check{e} \bullet ( \zeta^{-1} \theta ( f \bullet (\theta^{-1} \mu f')))
  =
  \zeta^{-1} \theta.((\theta^{-1} \zeta.\check{e})( f \bullet (\theta^{-1} \mu f')))
  =
  \zeta^{-1} \theta.(((\theta^{-1} \zeta.\check{e})f) (\theta^{-1} \mu.( (\mu^{-1} \theta. f) f')))\\
  = \
  &\zeta^{-1} \mu.( \mu^{-1} \theta.( 
  ((\theta^{-1} \zeta.\check{e})f) (\theta^{-1} \mu.( (\mu^{-1} \theta. f) f')))
  =
  \zeta^{-1} \mu.( \mu^{-1} \theta. ((\theta^{-1} \zeta.\check{e})f) (\mu^{-1} \theta. f) f')\\
  = \
  &\zeta^{-1} \mu.( \mu^{-1} \theta. ((\theta^{-1} \zeta.\check{e})f) f')
  =
  \zeta^{-1} \mu.( \mu^{-1} \zeta.( \zeta^{-1} \theta. ((\theta^{-1} \zeta.\check{e})f)) f')
  =
  \zeta^{-1} \mu.( \mu^{-1} \zeta.( \check{e} \bullet (\zeta^{-1} \theta f)) f')\\
  = \
  &(\check{e} \bullet (\zeta^{-1} \theta f)) \bullet (\zeta^{-1} \mu f').
\end{align*}
\eproof
Now, to prove c), take $e_0 \bullet (\zeta^{-1} \eta_1 e_1) \bullet \dotso \bullet (\zeta^{-1} \eta_{l-1} e_{l-1}) \in E_{\zeta, \theta}$ and $f_0 \bullet (\theta^{-1} \mu_1 f_1) \bullet \dotso \bullet (\theta^{-1} \mu_{n-1} f_{n-1}) \in E_{\theta,\nu}$. Then we have
\begin{align*}
  &(e_0 \bullet (\zeta^{-1} \eta_1 e_1) \bullet \dotso \bullet (\zeta^{-1} \eta_{l-1} e_{l-1}))
  \bullet
  (\zeta^{-1} \theta
  (f_0 \bullet (\theta^{-1} \mu_1 f_1) \bullet \dotso \bullet (\theta^{-1} \mu_{n-1} f_{n-1}))\\
  = \
  &e_0 \bullet (\zeta^{-1} \eta_1 e_1) \bullet \dotso \bullet (\zeta^{-1} \eta_{l-1} e_{l-1})
  \bullet
   (\zeta^{-1} \theta
  (f_0 \bullet (\theta^{-1} \mu_1 f_1) \bullet \dotso \bullet (\theta^{-1} \mu_{n-2} f_{n-2}))
  \bullet
  (\zeta^{-1} \mu_{n-1} f_{n-1})\\
  = \
  &\dotso
  = 
  e_0 \bullet (\zeta^{-1} \eta_1 e_1) \bullet \dotso \bullet (\zeta^{-1} \eta_{l-1} e_{l-1}) 
  \bullet
  (\zeta^{-1} \theta f_0) \bullet (\zeta^{-1} \mu_1 f_1) \bullet \dotso \bullet (\zeta^{-1} \mu_{n-1} f_{n-1})   
\end{align*}
which lies in $E_{\zeta,\nu}$, as desired. Here we used Lemma~\ref{lem:eff'}.
\setlength{\parindent}{0cm} \setlength{\parskip}{0.5cm}

Finally, let us prove b). Take $e_0 \bullet (\zeta^{-1} \eta_1 e_1) \bullet \dotso \bullet (\zeta^{-1} \eta_{l-1} e_{l-1}) \in E_{\zeta, \theta}$, and proceeding inductively on $l$ (the case $l=1$ being covered by property b) for $\cE_{\zeta,\theta}$), we have
\begin{align}
  &\theta^{-1} \zeta.(e_0 \bullet (\zeta^{-1} \eta_1 e_1) \bullet \dotso \bullet (\zeta^{-1} \eta_{l-1} e_{l-1}))
  =
  \theta^{-1} \zeta.( \zeta^{-1} \eta_{l-1}. (e_{l-1} ( \eta_{l-1}^{-1} \zeta.( e_0 \bullet (\zeta^{-1} \eta_1 e_1) \bullet \dotso)))) \nonumber \\
  = \
  &\theta^{-1} \eta_{l-1}. (e_{l-1} ( \eta_{l-1}^{-1} \zeta.( e_0 \bullet (\zeta^{-1} \eta_1 e_1) \bullet \dotso)))
  =
  \theta^{-1} \eta_{l-1}. (\eta_{l-1}^{-1} \theta. \theta^{-1} \eta_{l-1}. e_{l-1} ( \eta_{l-1}^{-1} \zeta.( e_0 \bullet (\zeta^{-1} \eta_1 e_1) \bullet \dotso))) \nonumber \\
\label{e:FinalProduct_b}
  = \
  &(\theta^{-1} \eta_{l-1}. e_{l-1}) \bullet (\theta^{-1} \eta_{l-1}( \eta_{l-1}^{-1} \zeta.( e_0 \bullet (\zeta^{-1} \eta_1 e_1) \bullet \dotso))).    
\end{align}
By induction hypothesis, $\eta_{l-1}^{-1} \zeta.( e_0 \bullet (\zeta^{-1} \eta_1 e_1) \bullet \dotso)$ lies in $E_{\eta_{l-1},\zeta}$, so that the term in \eqref{e:FinalProduct_b} lies in $E_{\theta,\zeta}$ by property c), which we have already established. Hence we have shown $\theta^{-1} \zeta. E_{\zeta,\theta} \subseteq E_{\theta,\zeta}$, as desired. This completes the proof of Lemma~\ref{lem:cE-E}.
\eproof

\bproof[Proof of Lemma~\ref{lem:E-FinDimAlg}]
It is clear that a) implies that \eqref{e:AlgE} defines an $F$-invariant subspace. It is *-invariant because
$$
  (d \delta_{\zeta^{-1}\eta} \otimes \ve_{\zeta,\eta})^*
  = (\delta_{\eta^{-1}\zeta} d \delta_{\zeta^{-1}\eta}) \delta_{\eta^{-1}\zeta} \otimes \ve_{\eta,\zeta} 
  = (\eta^{-1}\zeta.d) \, \delta_{\eta^{-1}\zeta} \otimes \ve_{\eta,\zeta},
$$
and $\eta^{-1}\zeta.d$ lies in $E_{\eta,\zeta}$ by b). To see that \eqref{e:AlgE} defines a subalgebra, we compute
\begin{align*}
  &(d \delta_{\zeta^{-1}\eta} \otimes \ve_{\zeta,\eta})
  (e \delta_{\eta^{-1}\theta} \otimes \ve_{\eta,\theta})
  = (d \delta_{\zeta^{-1}\eta} e \delta_{\eta^{-1}\theta}) \otimes \ve_{\zeta,\theta}
  = (d \delta_{\zeta^{-1}\eta} e \delta_{\eta^{-1}\zeta} \delta_{\zeta^{-1}\theta}) \otimes \ve_{\zeta,\theta}\\
  = \
  &(\delta_{\zeta^{-1}\eta} (\delta_{\eta^{-1}\zeta} d \delta_{\zeta^{-1}\eta}) e \delta_{\eta^{-1}\zeta}) \delta_{\zeta^{-1}\theta}) \otimes \ve_{\zeta,\theta}
  =
  ((d \bullet (\zeta^{-1}\eta e)) \delta_{\zeta^{-1}\theta}) \otimes \ve_{\zeta,\theta}.  
\end{align*}
Here we used properties (1), (2) and (3) of covariant representations of partial dynamical systems from \cite[Definition of covariant representations in \S~2]{McCl}. By c), we know that $d \bullet (\zeta^{-1}\eta e)$ lies in $E_{\zeta,\theta}$. This shows that \eqref{e:AlgE} defines a multiplicatively closed subspace, as desired.
\eproof

With these lemmas, we are now ready to construct inductive limit decompositions of $A$ and $\cA$. Recall that $F$ is a fixed finite subgroup of $G$. Let $\Sigma_i$, $i = 1, 2, 3, \dotsc$, be an increasing family of finite $F$-invariant subsets of $G$ such that $G = \bigcup_i \Sigma_i$. For each $i$, choose a family $\big\{\fE^{(i)}_{\zeta,\eta}\big\}_{\zeta,\eta \in \Sigma_i}$ of finite subsets $\fE^{(i)}_{\zeta,\eta}$ of $E_{\zeta^{-1}\eta}$ which is increasing in $i$ such that $E_{\zeta^{-1}\eta} = \bigcup_i \fE^{(i)}_{\zeta,\eta}$ for all $\zeta, \eta \in G$. First follow the procedure described right before Lemma~\ref{lem:fE-cE} to construct families $\big\{\cE^{(i)}_{\zeta,\eta}\big\}_{\zeta,\eta \in \Sigma}$ of finite subsemilattices $\cE^{(i)}_{\zeta,\eta}$ of $E_{\zeta^{-1}\eta}$ containing $\fE^{(i)}_{\zeta,\eta}$. Then follow the procedure described right before Lemma~\ref{lem:cE-E} to construct families $\big\{E^{(i)}_{\zeta,\eta}\big\}_{\zeta,\eta \in \Sigma}$ of finite subsemilattices $E^{(i)}_{\zeta,\eta}$ of $E_{\zeta^{-1}\eta}$ containing $\cE^{(i)}_{\zeta,\eta}$. Finally, set
$$
  A_i \defeq \lspan(\big\{d \delta_{\zeta^{-1}\eta} \otimes \ve_{\zeta,\eta}: \: d \in E^{(i)}_{\zeta,\eta}, \, \zeta, \eta \in \Sigma_i\big\}).
$$
By Lemma~\ref{lem:E-FinDimAlg}, $A_i$ are finite-dimensional sub-C*-algebras of $A$. By construction, $A_i$ form an increasing sequence such that $A = \ilim_i A_i$, or more precisely $A = \overline{\bigcup_i A_i}$. Moreover, define
$$
  \cA_i \defeq \lspan(\big\{\ve_{(d,\zeta),(\eta^{-1}\zeta.d,\eta)}: \: \zeta, \eta \in \Sigma_i, \, d \in E^{(i)}_{\zeta,\eta}\big\}).
$$
By construction, $\cA_i$ are finite-dimensional sub-C*-algebras of $\cA$, which form an increasing sequence such that $\cA = \ilim_i \cA_i$, or more precisely $\cA = \overline{\bigcup_i \cA_i}$.

What is more, these inductive limit decompositions are $F$-equivariant, so that we obtain $A \rtimes F = \ilim_i A_i \rtimes F$ and $\cA \rtimes F = \ilim_i \cA_i \rtimes F$. Since $A_i \rtimes F$ and $\cA_i \rtimes F$ are finite-dimensional for all $i$, we obtain
\blemma
\label{lem:AxF=AF}
For every finite subgroup $F$ of $G$, $A \rtimes F$ and $\cA \rtimes F$ are AF-algebras.
\elemma

Next, we define
$$
  \Phi_i: \: \cA_i \to \cK(\ell^2E) \otimes A_i, \, \ve_{(d,\zeta), \, (\eta^{-1}\zeta.d,\eta)} \ma \ve_{d, \, \eta^{-1}\zeta.d} \otimes (d \delta_{\zeta^{-1}\eta} \otimes \ve_{\zeta,\eta}).
$$
It is clear that $\Phi_i$ is $F$-equivariant, so that $\bm{x}_i \defeq [\Phi_i]$ defines an element in $KK^F(\cA_i,\cK(\ell^2E) \otimes A_i)$. By construction, we have a commutative diagram
\begin{equation*}
\begin{tikzcd}
\cA_i \ar["\Phi_i"]{r} \ar[hookrightarrow]{d} & \cK(\ell^2E) \otimes A_i \ar[hookrightarrow]{d} \\
\cA \ar["\Phi"]{r} & \cK(\ell^2E) \otimes A
\end{tikzcd}
\end{equation*}
where the vertical arrows are the canonical inclusions. Forming crossed products by $F$, we obtain the commutative diagram
\begin{equation*}
\begin{tikzcd}
\cA_i \rtimes F \ar["\Phi_i \rtimes F"]{r} \ar[hookrightarrow]{d} & \cK(\ell^2E) \otimes (A_i \rtimes F) \ar[hookrightarrow]{d} \\
\cA \rtimes F \ar["\Phi \rtimes F"]{r} & \cK(\ell^2E) \otimes (A \rtimes F)
\end{tikzcd}
\end{equation*}
Here we used the observation that $F$ acts trivially on the tensor factor $\cK(\ell^2E)$, so that we can pull it out. So we see that $\Phi \rtimes F = \ilim_i \Phi_i \rtimes F$. By continuity of K-theory, this means that if $\Phi_i \rtimes F$ induce isomorphisms in K-theory, then so does $\Phi \rtimes F$, i.e., taking Kasparov product with $j_F(\res^G_F(\bm{x}))$ is an isomorphism. So in summary, we obtain
\bprop
\label{prop:xi-->x}
If $\bm{x}_i = [\Phi_i]$ is invertible in $KK^F(\cA_i,\cK(\ell^2E) \otimes A_i)$ for all $i$, then taking Kasparov product with $j_F(\res^G_F(\bm{x}))$ induces an isomorphism $\sqcup \otimes j_F(\res^G_F(\bm{x})): \: K_*(\cA \rtimes F) \cong K_*((\cK(\ell^2E) \otimes A) \rtimes F)$.
\eprop

\subsection{KK-equivalences for finite subgroups}
\label{ss:KKEqFiniteSubg}

In the previous subsection, we have explained why it suffices to show that the $KK^F$-elements $\bm{x}_i = [\Phi_i]$ are invertible in $KK^F(\cA_i,\cK(\ell^2E) \otimes A_i)$. Our goal now is to show precisely this.

We start with the following observation, which is straightforward to check.
\blemma
\label{lem:cA=A}
We have an $F$-equivariant isomorphism 
$$
  \Psi_i: \: \cA_i \to A_i, \, \ve_{(d,\zeta),(\eta^{-1}\zeta.d,\eta)} \ma \Big(d - \bigvee_{e \in E^{(i)}_{\zeta,\eta}, \, e \lneq d} e\Big) \, \delta_{\zeta^{-1} \eta} \otimes \ve_{\zeta,\eta}.
$$
Its inverse is given by
$$
  A_i \to \cA_i, \, d \delta_{\zeta^{-1}\eta} \otimes \ve_{\zeta,\eta} \ma \sum_{e \in E^{(i)}_{\zeta,\eta}, \, e \leq d} \ve_{(e,\zeta),(\eta^{-1}\zeta.e,\eta)}.
$$
\elemma

Therefore, to show that $\bm{x}_i = [\Phi_i]$ is invertible it suffices to show that the following composite
\begin{equation*}
\begin{tikzcd}
\cA_i \ar["\Phi_i"]{r} & \cK(\ell^2E) \otimes A_i \ar["\id \otimes \Psi_i"]{r} & \cK(\ell^2E) \otimes \cA_i \\
\ve_{(d,\zeta),(\eta^{-1}\zeta.d,\eta)} \ar[maps to]{r} & \ve_{d,\eta^{-1}\zeta.d} \otimes (d \delta_{\zeta^{-1}\eta} \otimes \ve_{\zeta,\eta}) \ar[maps to]{r} & \ve_{d,\eta^{-1}\zeta.d} \otimes \sum_{e \in E^{(i)}_{\zeta,\eta}, \, e \leq d} \ve_{(e,\zeta),(\eta^{-1}\zeta.e,\eta)}
\end{tikzcd}
\end{equation*}
yields a $KK^F$-equivalence. Consider
\begin{eqnarray*}
  &&I: \: \cA_i \to \cK(\ell^2E) \otimes \cA_i, \, \ve_{(d,\zeta),(\eta^{-1}\zeta.d,\eta)} \ma \ve_{d,\eta^{-1}\zeta.d} \otimes \ve_{(d,\zeta),(\eta^{-1}\zeta.d,\eta)}\\
  &&\rho: \: \cA_i \to \cK(\ell^2E) \otimes \cA_i, \, \ve_{(d,\zeta),(\eta^{-1}\zeta.d,\eta)} \ma \ve_{d,\eta^{-1}\zeta.d} \otimes \sum_{e \in E^{(i)}_{\zeta,\eta}, \, e \lneq d} \ve_{(e,\zeta),(\eta^{-1}\zeta.e,\eta)}
\end{eqnarray*}
$I$ and $\rho$ are $F$-equivariant homomorphisms, and they are orthogonal because $\ve_{(d,\zeta),(\eta^{-1}\zeta.d,\eta)} \, \ve_{(e,\zeta),(\eta^{-1}\zeta.e,\eta)} = 0$ for all $e \lneq d$. Moreover, it is clear that $(id \otimes \Psi_i) \circ \Phi_i = I + \rho$.

Fix $f \in E\reg$. Then $\ve_{f,f} \otimes \id: \: \cA_i \to \cK(\ell^2E) \otimes \cA_i, \, a \ma \ve_{f,f} \otimes a$ is an invertible element in $KK^F(\cA_i,\cK(\ell^2E) \otimes \cA_i)$. This gives us a way to identify $KK^F(\cA_i,\cK(\ell^2E) \otimes \cA_i)$ with $KK^F(\cA_i,\cA_i)$.
\blemma
\label{lem:I=1}
Upon identifying $KK^F(\cA_i,\cK(\ell^2E) \otimes \cA_i)$ with $KK^F(\cA_i,\cA_i)$ as above, we have that $[I] = 1$ in $KK^F$.
\elemma
\setlength{\parindent}{0cm} \setlength{\parskip}{0cm}

\bproof
Consider
$$
  V \defeq \sum_{\zeta \in \Sigma_i, \, f \in E^{(i)}_{\zeta,\zeta}} \ve_{f,f} \otimes \ve_{(f,\zeta),(f,\zeta)} \ 
  + \
  \sum_{\zeta \in \Sigma_i, \, d \in E^{(i)}_{\zeta,\zeta}, \, d \neq f} \ve_{f,d} \otimes \ve_{(d,\zeta),(d,\zeta)} \ 
  + \
  \sum_{\zeta \in \Sigma_i, \, d \in E^{(i)}_{\zeta,\zeta}, \, d \neq f} \ve_{d,f} \otimes \ve_{(d,\zeta),(d,\zeta)}.
$$
$V$ is a self-adjoint partial isometry. Set $U \defeq W + (1-W^2) \in \cU(\cM(\cK(\ell^2E) \otimes \cA_i))$. Then $U$ is a self-adjoint unitary. We claim that $U I U = \ve_{f,f} \otimes \id_{\cA_i}$. Indeed, we have
\begin{align*}
  &U I(\ve_{(d,\zeta),(\eta^{-1}\zeta.d,\eta)}) U 
  = U (\ve_{d,\eta^{-1}\zeta.d} \otimes \ve_{(d,\zeta),(\eta^{-1}\zeta.d,\eta)}) U\\
  = \
  &U (\ve_{d,\eta^{-1}\zeta.d} \otimes \ve_{(d,\zeta),(\eta^{-1}\zeta.d,\eta)}) (\ve_{\eta^{-1}\zeta.d,f} \otimes \ve_{(\eta^{-1}\zeta.d,\eta),(\eta^{-1}\zeta.d,\eta)})
  = U (\ve_{d,f} \otimes \ve_{(d,\zeta),(\eta^{-1}\zeta.d,\eta)}) \\
  = \
  &(\ve_{f,d} \otimes \ve_{(d,\zeta),(d,\zeta)})
  (\ve_{d,f} \otimes \ve_{(d,\zeta),(\eta^{-1}\zeta.d,\eta)})
  = \ve_{f,f} \otimes \ve_{(d,\zeta),(\eta^{-1}\zeta.d,\eta)}
  = (\ve_{f,f} \otimes \id)(\ve_{(d,\zeta),(\eta^{-1}\zeta.d,\eta)}).
\end{align*}
\eproof
\setlength{\parindent}{0cm} \setlength{\parskip}{0.5cm}

Let $L$ be the longest proper chain $d_1 \gneq d_2 \gneq \dotso \gneq d_{L-1} \gneq d_L$ in $\bigcup_{\zeta,\eta \in \Sigma_i} E^{(i)}_{\zeta,\eta}$.
\blemma
\label{lem:rhonilpotent}
Upon identifying $KK^F(\cA_i,\cK(\ell^2E) \otimes \cA_i)$ with $KK^F(\cA_i,\cA_i)$ as above, we have that $[\rho]^L = 0$ in $KK^F$.
\elemma
\setlength{\parindent}{0cm} \setlength{\parskip}{0cm}

\bproof
It suffices to show that
$$
  (\id^{\otimes (L-1)} \otimes \rho) \circ \dotso \circ (\id^{\otimes 2} \otimes \rho) \circ (\id \otimes \rho) \circ \rho = 0
$$
as a homomorphism $\cA_i \to \cK(\ell^2E)^{\otimes L} \otimes \cA_i$. The reason is that, writing $\cK \defeq \cK(\ell^2E)$ and $\ve \defeq \ve_{f,f}$, we have the following commutative diagram
\begin{equation*}
\begin{tikzcd}
\cA_i \ar["\rho"]{r} & \cK \otimes A_i \ar["\id \otimes \rho"]{r} & \cK \otimes \cK \otimes \cA_i \ar["\id \otimes \id \otimes \rho"]{r} & \dotso \ar["\id^{\otimes (L-1)} \otimes \rho"]{r} & \cK^{\otimes L} \otimes \cA_i\\
& \cA_i \ar["\rho"]{r} \ar["\ve \otimes \id"]{u} & \cK \otimes A_i \ar["\id \otimes \rho"]{r} \ar["\ve \otimes \id"]{u} & \dotso \ar["\id^{\otimes (L-2)} \otimes \rho"]{r} & \cK^{\otimes (L-1)} \otimes \cA_i \ar["\ve \otimes \id"]{u}\\
& & \dotso \ar{r} \ar{u} & \dotso \ar{r} & \dotso \ar{u}\\
& & & \cA_i \ar["\rho"]{r} \ar{u} & \cK \otimes \cA_i \ar{u}\\
& & & & \cA_i \ar["\ve \otimes \id"]{u}
\end{tikzcd}
\end{equation*}
It shows that $[\rho]^L$ and the $KK^F$-class given by $(\id^{\otimes (L-1)} \otimes \rho) \circ \dotso \circ (\id^{\otimes 2} \otimes \rho) \circ (\id \otimes \rho) \circ \rho$ differ only by a $KK^F$-equivalence. 
\setlength{\parindent}{0cm} \setlength{\parskip}{0.5cm}

We compute
\begin{eqnarray*}
  && \rho(\ve_{(d_1,\zeta),(\eta^{-1}\zeta.d_1,\eta)}) = \ve_{d_1,\eta^{-1}\zeta.d_1} \otimes \sum_{d_2 \lneq d_1} \ve_{(d_2,\zeta),(\eta^{-1}\zeta.d_2,\eta)}\\
  && ((\id \otimes \rho) \circ \rho)(\ve_{(d_1,\zeta),(\eta^{-1}\zeta.d_1,\eta)}) = \sum_{d_3 \lneq d_2 \lneq d_1} \ve_{d_1,\eta^{-1}\zeta.d_1} \otimes \ve_{(d_2,\zeta),(\eta^{-1}\zeta.d_2,\eta)} \otimes \ve_{(d_3,\zeta),(\eta^{-1}\zeta.d_3,\eta)}\\
  && \dotso\\\
  && ((\id^{\otimes (L-1)} \otimes \rho) \circ \dotso \circ \rho)(\ve_{(d_1,\zeta),(\eta^{-1}\zeta.d_1,\eta)})\\
  && = \sum_{d_{L+1} \lneq d_L \lneq \dotso \lneq d_2 \lneq d_1} \ve_{d_1,\eta^{-1}\zeta.d_1} \otimes \ve_{(d_2,\zeta),(\eta^{-1}\zeta.d_2,\eta)} \otimes \dotso \otimes \ve_{(d_L,\zeta),(\eta^{-1}\zeta.d_L,\eta)} \otimes \ve_{(d_{L+1},\zeta),(\eta^{-1}\zeta.d_{L+1},\eta)}\\
  && = 0
\end{eqnarray*}
because there are no proper chains of length greater than $L$.
\eproof
\setlength{\parindent}{0cm} \setlength{\parskip}{0.5cm}

\bcor
\label{cor:xiInv}
For all $i$, the element $\bm{x}_i = [\Phi_i]$ is invertible in $KK^F(\cA_i,\cK(\ell^2E) \otimes A_i)$.
\ecor
\setlength{\parindent}{0cm} \setlength{\parskip}{0cm}

\bproof
As explained above, it suffices to show that $[\id \otimes \Psi_i] \otimes [\Phi_i]$ is invertible in $KK^F(\cA_i,\cK(\ell^2E) \otimes A_i)$. We have $[\id \otimes \Psi_i] \otimes [\Phi_i] = [I] + [\rho]$ which by Lemma~\ref{lem:I=1} is given by $1 + [\rho]$ (upon identifying $KK^F(\cA_i,\cK(\ell^2E) \otimes \cA_i)$ with $KK^F(\cA_i,\cA_i)$), which then by Lemma~\ref{lem:rhonilpotent} is invertible with inverse $\sum_{l=0}^{L-1} (-1)^l [\rho]^l$.
\eproof
\setlength{\parindent}{0cm} \setlength{\parskip}{0.5cm}

\subsection{Homomorphisms inducing KK-equivalence}
\label{ss:HomIndKKEq}

Let us first analyse $\cA \rtimes_r G$. We first need some notation. Given $d \in E\reg$, let $G(d) \defeq \menge{\gamma \in G}{d \in E_{\gamma^{-1}}}$ and $G_d \defeq \menge{\gamma \in G(d)}{\gamma.d = d}$. Choose a set of representatives $\fD$ for $G \backslash E\reg$, i.e., a subset $\fD \subseteq E\reg$ such that for every $e \in E\reg$, there exists a unique $d \in \fD$ and some $\gamma \in G$ with $e = \gamma.d$. For $d \in \fD$, define a subalgebra $\cA_d$ of $\cA$ by $\cA_d \defeq \clspan(\menge{\ve_{(e,\zeta),(\eta^{-1}\zeta.e,\eta)}}{\zeta, \eta \in G, \, e \in G(d).d \cap E_{\zeta^{-1}\eta}})$. We view $\cA_d$ as a subalgebra of $\cK(\ell^2(G(d).d \times G))$. Clearly, we have $\cA = \bigoplus_{d \in \fD} \cA_d$. Moreover, recall that $G$ acts on $\cA$ by $\Ad(1 \otimes \lambda)$. By construction, $\cA_d$ is $G$-invariant for all $d \in \fD$. Thus we will focus on one summand $\cA_d$ for a fixed $d \in \fD$.

Let $\fR \subseteq G$ be a set of representatives for $G / G_d$, and let $\fr: \: G/G_d \to \fR$ be a split for the canonical map $\fR \into G \onto G/G_d$ (i.e., we have $[\fr[\gamma]] = [\gamma]$ for all $\gamma \in G$). Now consider the map
$$
  G(d).d \times G \to G(d)/G_d \times G/G_d \times G_d, \, (\gamma.d,\zeta) \ma ([\gamma], [\zeta \gamma], (\fr[\zeta \gamma])^{-1} \zeta \fr[\gamma]).
$$
It is a bijection with inverse 
$$
  G(d)/G_d \times G/G_d \times G_d \to G(d).d \times G, \, ([\gamma], [\tau], \mu) \ma (\fr[\gamma].d, \fr[\tau] \mu \fr[\gamma]^{-1}).
$$
This bijection induces a unitary
$$
  \ell^2(G(d).d \times G) \cong \ell^2(G(d)/G_d \times G/G_d \times G_d) \cong \ell^2(G(d)/G_d) \otimes \ell^2(G/G_d) \otimes \ell^2G_d.
$$
Conjugation by this unitary yields an isomorphism
$$
  \Xi_d: \: \cA_d \cong \cK(\ell^2(G(d)/G_d)) \otimes C_0(G/G_d) \otimes \cK(\ell^2G_d)
$$
with
$$
  \Xi_d(\ve_{(e,\zeta),(\eta^{-1}\zeta.e,\eta)}) = \ve_{[\gamma],[\eta^{-1}\zeta \gamma]} \otimes \ve_{[\zeta \gamma], [\zeta \gamma]} \otimes \ve_{(\fr[\zeta \gamma])^{-1} \zeta \fr[\gamma], (\fr[\zeta \gamma])^{-1} \eta \fr[\eta^{-1}\zeta \gamma]}
$$
where $e = \gamma.d$. 

Now let $l$ denote the $G$-action on $C_0(G/G_d)$ given by left translation.
\blemma
\label{lem:Ad1lambdaCC1l1}
The $G$-actions $\Ad(1 \otimes \lambda): \: G \curvearrowright \cA_d$ and $\id \otimes l \otimes \id: \: G \curvearrowright \cK(\ell^2(G(d)/G_d)) \otimes C_0(G/G_d) \otimes \cK(\ell^2G_d)$ are cocycle conjugate.
\elemma
\setlength{\parindent}{0cm} \setlength{\parskip}{0cm}

\bproof
For $g \in G$, let $w_g$ be the unitary $\ell^2(G(d)/G_d) \otimes \ell^2(G/G_d) \otimes \ell^2G_d \cong \ell^2(G(d)/G_d) \otimes \ell^2(G/G_d) \otimes \ell^2G_d$ induced by the bijection
$$
  G(d)/G_d \times G/G_d \times G_d \to G(d)/G_d \times G/G_d \times G_d, \, ([\gamma],[\tau],\mu) \ma ([\gamma],[\tau],(\fr[\tau])^{-1} g \fr[g^{-1} \tau] \mu).
$$
Let us show that
$$
  \Xi_d \circ \Ad(1 \otimes \lambda_g) \circ \Xi_d^{-1} = w_g (\id \otimes l_g \otimes \id) w_g^*.
$$
$\Xi_d \circ \Ad(1 \otimes \lambda_g) \circ \Xi_d^{-1}$ is given by conjugation with the unitary which corresponds to the following bijection:
\begin{equation*}
\adjustbox{scale=0.9333,center}{
\begin{tikzcd}
  G(d)/G_d \times G/G_d \times G_d \ar{r} & G(d).d \times G \ar{r} & G(d).d \times G \ar{r} & G(d)/G_d \times G/G_d \times G_d\\
  ([\gamma],[\tau],\mu) \ar[maps to]{r} & (\fr[\gamma].d,\fr[\tau] \mu (\fr[\gamma])^{-1}) \ar[maps to]{r} & (\fr[\gamma].d, g \fr[\tau] \mu (\fr[\gamma])^{-1}) \ar[maps to]{r} & ([\gamma], [g \tau], (\fr[g \tau])^{-1} g \fr[\tau] \mu)
\end{tikzcd}
}
\end{equation*}
$w_g (\id \otimes l_g \otimes \id) w_g^*$ is given by conjugation with the unitary which corresponds to the following bijection:
\begin{equation*}
\begin{tikzcd}
  G(d)/G_d \times G/G_d \times G_d \ar{r} & G(d)/G_d \times G/G_d \times G_d \ar{r} & G(d)/G_d \times G/G_d \times G_d\\
  ([\gamma],[\tau],\mu) \ar[maps to]{r} & ([\gamma],[g\tau],\mu) \ar[maps to]{r} & ([\gamma], [g \tau], (\fr[g \tau])^{-1} g \fr[\tau] \mu)
\end{tikzcd}
\end{equation*}
This shows $\Xi_d \circ \Ad(1 \otimes \lambda_g) \circ \Xi_d^{-1} = w_g (\id \otimes l_g \otimes \id) w_g^*$, as desired.
\setlength{\parindent}{0cm} \setlength{\parskip}{0.5cm}

It remains to show that $g \ma w_g$ is an ($\id \otimes l \otimes \id$)-cocycle. For that purpose, observe that $w_g (\id \otimes l_g \otimes \id)(w_h)$ is the unitary induced by the bijection
\begin{equation*}
\begin{tikzcd}
  G(d)/G_d \times G/G_d \times G_d \ar{r} & G(d)/G_d \times G/G_d \times G_d \ar{r} & G(d)/G_d \times G/G_d \times G_d\\
  ([\gamma],[\tau],\mu) \ar[maps to]{r} & ([\gamma],[\tau],(\fr[g^{-1} \tau])^{-1} h \fr[h^{-1} g^{-1} \tau] \mu) \ar[maps to]{r} & ([\gamma], [\tau], (\fr[\tau])^{-1} g h \fr[h^{-1} g^{-1} \tau] \mu)
\end{tikzcd}
\end{equation*}
which is precisely the bijection corresponding to $w_{gh}$. Hence $w_{gh} = w_g (\id \otimes l_g \otimes \id)(w_h)$, as desired.
\eproof
\setlength{\parindent}{0cm} \setlength{\parskip}{0.5cm}

Let us denote the canonical unitaries in (the multiplier algebra of) $\cA \rtimes_r G$ implementing the $G$-action by $\Delta_g$, $g \in G$. For $d \in \fD$, define the homomorphism $\kappa_d: \: C^*_{\lambda}(G_d) \to \cA \rtimes_r G, \, \lambda_g \ma \ve_{(d,1),(d,g)} \Delta_g$. 
\bprop
\label{prop:C=cAG}
$$
  \Big(\sum_{d \in \fD} \kappa_d\Big): \Big(\bigoplus_{d \in \fD} C^*_{\lambda}(G_d)\Big) \to \cA \rtimes_r G
$$
induces a Morita equivalence between $\bigoplus_{d \in \fD} C^*_{\lambda}(G_d)$ and $\cA \rtimes_r G$.
\eprop
\setlength{\parindent}{0cm} \setlength{\parskip}{0cm}

\bproof
Consider the composition
\begin{equation}
\label{e:CtoKCKtoA}
  \bigoplus_{d \in \fD} C^*_{\lambda}(G_d) 
  \to 
  \bigoplus_{d \in \fD} \big( \cK(\ell^2(G(d)/G_d)) \otimes C_0(G/G_d) \otimes \cK(\ell^2G_d) \big) \rtimes_{\id \otimes l \otimes \id, \, r} G
  \overset{\cong}{\lori}
  \bigoplus_{d \in \fD} \cA_d \rtimes_r G \cong \cA \rtimes_r G,  
\end{equation}
where the first map sends $\lambda_g \in C^*_{\lambda}(G_d)$ to $(\ve_{[1],[1]} \otimes \ve_{[1],[1]} \otimes \ve_{1,1}) \, \Delta_g$, the second map is given by the isomorphism sending $a \Delta_g \in \big( \cK(\ell^2(G(d)/G_d)) \otimes C_0(G/G_d) \otimes \cK(\ell^2G_d) \big) \rtimes_{\id \otimes l \otimes \id, \, r} G$ to $\Xi_d^{-1}(a) w_g^{-1} \Delta_g$, and the third map is the canonical isomorphism. Now it is clear that the composition in \eqref{e:CtoKCKtoA} induces a Morita equivalence between $\bigoplus_{d \in \fD} C^*_{\lambda}(G_d)$ and $\cA \rtimes_r G$ (this is as in \cite[\S~3]{CEL2}). In addition, it is easy to check that the composition in \eqref{e:CtoKCKtoA} coincides with $\sum_{d \in \fD} \kappa_d$.
\eproof
\setlength{\parindent}{0cm} \setlength{\parskip}{0.5cm}

If we now define the homomorphisms $\fk_d: \: C^*_{\lambda}(G_d) \to \cK(\ell^2E) \otimes \fA \rtimes_r G, \, \lambda_g \ma (d \delta_g \otimes \ve_{1,g}) \Delta_g$, then the following is easy to check:
\blemma
\label{lem:CtofAtoA}
We have a commutative diagram
\begin{equation*}
\begin{tikzcd}
  \bigoplus_{d \in \fD} C^*_{\lambda}(G_d) \ar["\sum_d \ve_{d,d} \otimes \, \fk_d"]{rr} \ar{d}[swap]{\sum_d \kappa_d} & & \cK(\ell^2E) \otimes \fA \rtimes_r G \ar[hookrightarrow]{d}\\
  \cA \rtimes_r G \ar["\Phi \rtimes_r G"]{rr} & & \cK(\ell^2E) \otimes A \rtimes_r G
\end{tikzcd}
\end{equation*}
where the right vertical arrow is given by the canonical inclusion.
\elemma

Given $d \in \fD$, define the homomorphism $\iota_d: \: C^*_{\lambda}(G_d) \to D \rtimes_r G, \, \lambda_g \ma d \delta_g$. Let $[M \rtimes_r G]$ be the element in $KK(D \rtimes_r G, \fA \rtimes_r G)$ induced by the $(D \rtimes_r G)$-$(\fA \rtimes_r G)$-imprimitivity bimodule $M \rtimes_r G$ introduced in \eqref{e:MxG} in \S~\ref{ss:MoritaEnv}. 
\bprop
\label{prop:iotaM=fk}
For every $d \in \fD$, we have $[\iota_d] \otimes [M \rtimes_r G] = [\fk_d]$ in $KK(C^*_{\lambda}(G_d),\fA \rtimes_r G)$.
\eprop
\setlength{\parindent}{0cm} \setlength{\parskip}{0cm}

\bproof
The Kasparov product $[\iota_d] \otimes [M \rtimes_r G]$ is represented by the right Hilbert $\fA \rtimes_r G$-module $D \rtimes_r G \otimes_{D \rtimes_r G} M \rtimes_r G \cong M \rtimes_r G$, with right Hilbert $\fA \rtimes_r G$-module as described in \S~\ref{ss:MoritaEnv} and left $C^*_{\lambda}(G_d)$-action given by $\lambda_g ((e \delta_{\eta} \otimes \ve_{1,\eta}) \Delta_h) = (g.(de) \delta_{g \eta} \otimes \ve_{1,g \eta}) \Delta_{gh}$. The KK-element $[\fk_d]$ is represented by $\fA \rtimes_r G$, viewed as a right Hilbert $\fA \rtimes_r G$ in the obvious way, with left $C^*_{\lambda}(G_d)$-action given by
$$
  \lambda_g \, ((e \delta_{\zeta^{-1}\eta} \otimes \ve_{\zeta,\eta}) \Delta_h) = (d \delta_g \otimes \ve_{1,g}) (e \delta_{\zeta^{-1}\eta} \otimes \ve_{g\zeta,g\eta}) \Delta_{gh}
  =
  \bfa
    0 & {\rm if} \ \zeta \neq 1,\\
    ((g.(de)) \delta_{g \eta} \otimes \ve_{1,g \eta} ) \Delta_{gh} & {\rm if} \ \zeta = 1.
  \efa
$$
Without changing the classes in $KK(C^*_{\lambda}(G_d),\fA \rtimes_r G)$, we can replace the Hilbert modules $M \rtimes_r G$ by $C^*_{\lambda}(G_d) \cdot (M \rtimes_r G)$ and $\fA \rtimes_r G$ by $C^*_{\lambda}(G_d) \cdot (\fA \rtimes_r G)$. Then it is straightforward to check that 
$$
  C^*_{\lambda}(G_d) \cdot (M \rtimes_r G) \to C^*_{\lambda}(G_d) \cdot (\fA \rtimes_r G), \, (e \delta_{\eta} \otimes \ve_{1,\eta}) \Delta_g \ma (e \delta_{\eta} \otimes \ve_{1,\eta}) \Delta_g
$$
gives rise to an isomorphism of right Hilbert $\fA \rtimes_r G$-modules which intertwines the left $C^*_{\lambda}(G_d)$-actions.
\eproof
\setlength{\parindent}{0cm} \setlength{\parskip}{0.5cm}

\subsection{Proofs of main theorems}
\label{ss:Pfs}

Let us now prove Theorems~\ref{thm:MainInvSgp}, \ref{thm:MainPartialCroPro} and Corollaries~\ref{cor:MainP}, \ref{cor:MainLCM}. We actually prove more precise statements. Let us keep the same notations as in previous subsections of \S~\ref{s:KMoritaEnv}. For $d \in \fD$, let $i_d: \: C^*_{\lambda}(S_d) \to C^*_{\lambda}(S)$ be the homomorphism induced by the canonical embedding $S_d \into S$.
\btheo
\label{thm:InvSgp}
$ $
\setlength{\parindent}{0cm} \setlength{\parskip}{0cm}

\begin{enumerate}
\item[(I)] If $G$ satisfies the Baum-Connes conjecture for $\cA$ and $A$, then $\sum_{d \in \fD} i_d$ induces a K-theory isomorphism $\bigoplus_{d \in \fD} K_*(C^*_{\lambda}(S_d)) \cong K_*(C^*_{\lambda}(S))$.
\item[(II)] If $G$ satisfies the strong Baum-Connes conjecture  in the sense of \cite[Definition~3.4.17]{CELY}, then $\sum_{d \in \fD} i_d$ induces a KK-equivalence between $\bigoplus_{d \in \fD} C^*_{\lambda}(S_d)$ and $C^*_{\lambda}(S)$.
\end{enumerate}
\etheo
\setlength{\parindent}{0cm} \setlength{\parskip}{0cm}

\bproof
For every $d \in \fD$, $\sigma$ induces an isomorphism $S_d \cong G_d$ and hence a C*-isomorphism $C^*_{\lambda}(S_d) \cong C^*_{\lambda}(G_d)$. Moreover, it is easy to check that this C*-isomorphism fits into the following commutative diagram:
\begin{equation*}
\begin{tikzcd}
  C^*_{\lambda}(S_d) \ar["i_d"]{r} \ar{d}[swap]{\cong} & C^*_{\lambda}(S)  \ar{d}{\cong}\\
  C^*_{\lambda}(G_d) \ar["\iota_d"]{r} & D \rtimes_r G
\end{tikzcd}
\end{equation*}
where the right vertical isomorphism is provided by \eqref{e:CS=DxG}. Hence it suffices to show that $\sum_{d \in \fD} \iota_d$ induces a K-theory isomorphism $\bigoplus_{d \in \fD} K_*(C^*_{\lambda}(G_d)) \cong K_*(D \rtimes_r G)$ in case (I) and a KK-equivalence between $\bigoplus_{d \in \fD} C^*_{\lambda}(G_d)$ and $D \rtimes_r G$ in case (II).
\setlength{\parindent}{0cm} \setlength{\parskip}{0.5cm}

Let $\bm{x} \defeq [\Phi] \in KK^G(\cA, \cK(\ell^2E) \otimes A)$ be as in \S~\ref{ss:DiscreteV}. Corollary~\ref{cor:xiInv} and Proposition~\ref{prop:xi-->x} imply that taking Kasparov product with $j_F(\res^G_F(\bm{x}))$ induces an isomorphism $\sqcup \otimes j_F(\res^G_F(\bm{x})): \: K_*(\cA \rtimes F) \cong K_*((\cK(\ell^2E) \otimes A) \rtimes F)$. Since $\cA \rtimes F$ and $(\cK(\ell^2E) \otimes A) \rtimes F \cong \cK(\ell^2E) \otimes (A \rtimes F)$ are AF-algebras by Lemma~\ref{lem:AxF=AF}, they satisfy the UCT, so that we actually obtain that $j_F(\res^G_F(\bm{x}))$ is a KK-equivalence between $\cA \rtimes F$ and $(\cK(\ell^2E) \otimes A) \rtimes F$.  Therefore, by the Going-Down principle as in Proposition~\ref{prop:GD}, we obtain that taking Kasparov product with $j_G(\bm{x})$ induces an isomorphism $\sqcup \otimes j_G(\bm{x}): \: K_*(\cA \rtimes_r G) \cong K_*((\cK(\ell^2E) \otimes A) \rtimes_r G)$ in case (I) and that $j_G(\bm{x})$ is a KK-equivalence between $\cA \rtimes_r G$ and $(\cK(\ell^2E) \otimes A) \rtimes_r G$ in case (II). Now Lemma~\ref{lem:CtofAtoA} and Proposition~\ref{prop:iotaM=fk} imply that the following diagram in KK commutes:
\begin{equation*}
\begin{tikzcd}
   & & \cK(\ell^2E) \otimes D \rtimes_r G \ar{d}{1 \otimes [M \rtimes_r G]}\\
  \bigoplus_{d \in \fD} C^*_{\lambda}(G_d) \ar{rr}[swap]{\big[ \sum_d \ve_{d,d} \otimes \, \fk_d \big]} \ar{d}[swap]{\big[ \sum_d \kappa_d \big]} \ar{urr}{\big[ \sum_d \ve_{d,d} \otimes \, \iota_d \big]} & & \cK(\ell^2E) \otimes \fA \rtimes_r G \ar{d}\\
  \cA \rtimes_r G \ar{rr}{j_G(\bm{x}) \, = \, [\Phi \rtimes_r G]} & & \cK(\ell^2E) \otimes A \rtimes_r G
\end{tikzcd}
\end{equation*}
where the lower right vertical arrow is given by the canonical inclusion. Since $\big[ \sum_d \kappa_d \big]$ is a KK-equivalence by Proposition~\ref{prop:C=cAG}, and because the two right vertical arrows are KK-equivalences by the discussion in \S~\ref{ss:MoritaEnv}, we deduce that $\sum_d \ve_{d,d} \otimes \, \iota_d$ must induce a K-theory isomorphism in case (I) and a KK-equivalence in case (II). Since for every $d \in \fD$, the KK-elements $[\iota_d]$ and $[\ve_{d,d} \otimes \, \iota_d \big]$ coincide upon identifying -- up to KK-equivalence -- $D \rtimes_r G$ with $\cK(\ell^2E) \otimes D \rtimes_r G$, we obtain the desired statement for $\sum_{d \in \fD} \iota_d$ and hence for $\sum_{d \in \fD} i_d$.
\eproof
Clearly, Theorem~\ref{thm:InvSgp} implies Theorem~\ref{thm:MainInvSgp}, once we observe that $G \backslash E\reg = S \backslash E\reg$.
\setlength{\parindent}{0cm} \setlength{\parskip}{0.5cm}

Now let us come to partial dynamical systems and their reduced crossed products. Let $G$ be a countable discrete group, $X$ a second countable totally disconnected locally compact Hausdorff space, and $G \curvearrowright X$ a partial dynamical system, given by $U_{g^{-1}} \to U_g, \, x \ma g.x$. Suppose that $\cV$ is a $G$-invariant regular basis for the compact open subsets of $X$. For $V \in \cV$, set $i_V: \: C^*_{\lambda}(G_V) \to C_0(X) \rtimes_r G, \, \lambda_g \ma \bm{1}_V \delta_g$. Here $G_V = \menge{g \in G}{g.V = V}$ and $\bm{1}_V$ denotes the characteristic function on $V$. Moreover, let $G \backslash \cV\reg$ denote the set of orbits under the $G$-action on the non-empty elements $\cV\reg$ of $\cV$. Apply construction ($**$) from \S~\ref{ss:PDS} to construct the inverse semigroup $S$, together with an idempotent pure partial homomorphism $\sigma: \: S\reg \to G$, attached to $G \curvearrowright X$ and $\cV$, and let $\cA$ and $A$ be the C*-algebras constructed in \S~\ref{ss:MoritaEnv} and \S~\ref{ss:DiscreteV}, respectively.
\btheo
\label{thm:PartialCroPro}
$ $
\setlength{\parindent}{0cm} \setlength{\parskip}{0cm}

\begin{enumerate}
\item[(I)] If $G$ satisfies the Baum-Connes conjecture for $\cA$ and $A$, then $\sum_{[V] \in G \backslash \cV\reg} i_V$ induces a K-theory isomorphism $\bigoplus_{[V] \in G \backslash \cV\reg} K_*(C^*_{\lambda}(G_V)) \cong K_*(C_0(X) \rtimes_r G)$.
\item[(II)] If $G$ satisfies the strong Baum-Connes conjecture in the sense of \cite[Definition~3.4.17]{CELY}, then $\sum_{[V] \in G \backslash \cV\reg} i_V$ induces a KK-equivalence between $\bigoplus_{[V] \in G \backslash \cV\reg} C^*_{\lambda}(G_V)$ and $C_0(X) \rtimes_r G$.
\end{enumerate}
\etheo
\setlength{\parindent}{0cm} \setlength{\parskip}{0cm}

\bproof
It is easy to check that for $S$ (with its semilattice $E$), $\sigma$ obtained by construction ($**$) from $G \curvearrowright X$ and $\cV$, we get a semilattice isomorphism $E \cong \cV$ respecting the partial $G$-actions. Now the theorem follows from the isomorphism in \eqref{e:CS=DxG} and Theorem~\ref{thm:InvSgp}.
\eproof
\setlength{\parindent}{0cm} \setlength{\parskip}{0.5cm}

Before we turn to semigroup C*-algebras, let us record a K-theory formula for reduced C*-algebras of left inverse hulls. Let $P$ be a subsemigroup of a countable group $G$. Let $I_l(P)$ be the left inverse hull of $P$, as introduced in \S~\ref{ss:InvSgp}. Let $\cA$ and $A$ be the C*-algebras constructed for $S = I_l(P)$ in \S~\ref{ss:MoritaEnv} and \S~\ref{ss:DiscreteV}. For $X \in \cJ_P\reg = E(I_l(P))\reg$, let $i_X: \: C^*_{\lambda}((I_l(P))_X) \to C^*_{\lambda}(I_l(P))$ be the homomorphism induced by the canonical embedding $(I_l(P))_X \into I_l(P)$.
\bcor
\label{cor:IlP}
$ $
\setlength{\parindent}{0cm} \setlength{\parskip}{0cm}

\begin{enumerate}
\item[(I)] If $G$ satisfies the Baum-Connes conjecture for $\cA$ and $A$, then $\sum_{[X] \in I_l(P) \backslash \cJ_P\reg} i_X$ induces a K-theory isomorphism $\bigoplus_{[X] \in I_l(P) \backslash \cJ_P\reg} K_*(C^*_{\lambda}((I_l(P))_X)) \cong K_*(C^*_{\lambda}(I_l(P)))$.
\item[(II)] If $G$ satisfies the strong Baum-Connes conjecture in the sense of \cite[Definition~3.4.17]{CELY}, then $\sum_{[X] \in I_l(P) \backslash \cJ_P\reg} i_X$ induces a KK-equivalence between $\bigoplus_{[X] \in I_l(P) \backslash \cJ_P\reg} C^*_{\lambda}((I_l(P))_X)$ and $C^*_{\lambda}(I_l(P))$.
\end{enumerate}
\ecor
\setlength{\parindent}{0cm} \setlength{\parskip}{0cm}

\bproof
Just apply Theorem~\ref{thm:InvSgp} to $S = I_l(P)$.
\eproof
\setlength{\parindent}{0cm} \setlength{\parskip}{0.5cm}

The case of semigroup C*-algebras can now be treated as a special case. Let $P$ be a left-cancellative semigroup, $\cJ_P\reg$ the set of non-empty constructible right ideals of $P$, as introduced in \S~\ref{ss:InvSgp}, $P \backslash \cJ_P\reg$ the set of equivalence classes of the equivalence relation on $\cJ_P\reg$ generated by $X \sim pX = \menge{px}{x \in X}$ for all $X \in \cJ_P\reg$ and $p \in P$, and, for $X \in \cJ_P\reg$, let $P_X$ be the group of bijections $X \to X$ which can be expressed as compositions of finitely many maps, each of which given by left multiplication by a fixed semigroup element or the set-theoretical inverse of such a left multiplication map. For $X \in \cJ_P\reg$, set $i_X: \: C^*_{\lambda}(P_X) \to C^*_{\lambda}(P), \, \lambda_g \ma V_g$. As in \S~\ref{ss:MoritaEnv} and \S~\ref{ss:DiscreteV}, construct C*-algebras $\cA$ and $A$ for the left inverse hull $S \defeq I_l(P)$ of $P$ (see \S~\ref{ss:InvSgp}).
\bcor
\label{cor:P}
Suppose that $P$ satisfies the independence condition from Definition~\ref{def:Ind}.
\setlength{\parindent}{0cm} \setlength{\parskip}{0cm}

\begin{enumerate}
\item[(I)] If $P$ embeds into a countable group which satisfies the Baum-Connes conjecture for $\cA$ and $A$, then $\sum_{[X] \in P \backslash \cJ_P\reg} i_X$ induces a K-theory isomorphism $\bigoplus_{[X] \in P \backslash \cJ_P\reg} K_*(C^*_{\lambda}(P_X)) \cong K_*(C^*_{\lambda}(P))$.
\item[(II)] If $P$ embeds into a countable group which satisfies the strong Baum-Connes conjecture in the sense of \cite[Definition~3.4.17]{CELY}, then $\sum_{[X] \in P \backslash \cJ_P\reg} i_X$ induces a KK-equivalence between $\bigoplus_{[X] \in P \backslash \cJ_P\reg} C^*_{\lambda}(P_X)$ and $C^*_{\lambda}(P)$.
\end{enumerate}
\ecor
\setlength{\parindent}{0cm} \setlength{\parskip}{0cm}

\bproof
This follows from Proposition~\ref{prop:P=S} and Corollary~\ref{cor:IlP}, once we observe that $P \backslash \cJ_P\reg = S \backslash E\reg$ and $P_X = (I_l(P))_X = S_X$ for all $X \in \cJ_P\reg$.
\eproof
\setlength{\parindent}{0cm} \setlength{\parskip}{0.5cm}

Let us further specialize to the case of right LCM monoids, i.e., monoids $P$ for which $\cJ_P\reg = \menge{pP}{p \in P}$. Let us keep the same notations as in Corollary~\ref{cor:P}, and denote by $P^*$ the group of invertible elements in $P$.
\bcor
\label{cor:LCM}
Let $P$ be a right LCM monoid.
\setlength{\parindent}{0cm} \setlength{\parskip}{0cm}

\begin{enumerate}
\item[(I)] If $P$ embeds into a countable group which satisfies the Baum-Connes conjecture for $\cA$ and $A$, then $i_P$ induces a K-theory isomorphism $K_*(C^*_{\lambda}(P^*)) \cong K_*(C^*_{\lambda}(P))$.
\item[(II)] If $P$ embeds into a countable group which satisfies the strong Baum-Connes conjecture in the sense of \cite[Definition~3.4.17]{CELY}, then $i_P$ induces a KK-equivalence between $C^*_{\lambda}(P^*)$ and $C^*_{\lambda}(P)$.
\end{enumerate}
\ecor
\setlength{\parindent}{0cm} \setlength{\parskip}{0cm}

\bproof
Left-cancellative right LCM monoids satisfy the independence condition by \cite[Lemma~5.6.31]{CELY}. Hence Corollary~\ref{cor:P} applies and we obtain the desired statement once we identify $P_P$ in the notation of Corollary~\ref{cor:P} with $P^*$. If we identify elements $u \in P^*$ with the left multiplication maps $P \to P, \, x \ma ux$, then it is clear that $P^* \subseteq P_P$. To prove the reverse inclusion, take $s \in P_P$. Let $1$ denote the identity element of $P$, and let $u \defeq s(1)$. Then, by \cite[Equation~(5.11)]{CELY}, we must have that $s(x) = ux$ for all $x \in P$. Thus $u$ must lie in $P^*$, and we obtain $P_P \subseteq P^*$, as desired. 
\eproof
\setlength{\parindent}{0cm} \setlength{\parskip}{0.5cm}

\section{Applications}
\label{s:App}

We start with applications of Corollary~\ref{cor:LCM} to semigroup C*-algebras of Artin monoids, Baumslag-Solitar monoids and one-relator monoids. We then discuss C*-algebras generated by right regular representations of $ax+b$-type semigroups of number-theoretic origin. Finally, we compute K-theory for reduced C*-algebras of inverse semigroups arising in the context of tilings.

\subsection{Semigroup C*-algebras of Artin monoids}
\label{ss:Artin}

First of all, let us discuss Artin monoids (sometimes also called Artin-Tits monoids). Let $\cS$ be a countable set. For every $a, b \in \cS$, let $m_{a,b} \in \gekl{2, 3, \dotsc} \cup \gekl{\infty}$ such that $m_{a,b} = m_{b,a}$. The Artin group $A_M$ of $M = (m_{a,b})_{a,b \in \cS}$ is the group given by the presentation 
$$
  A_M \defeq \spkl{\cS \ \vert \ \spkl{ab}^{m_{a,b}} = \spkl{ba}^{m_{b,a}} \ \forall \ a, b \in \cS}.
$$
Here $\spkl{ab}^{m_{a,b}}$ denotes the alternating word $abab \dotsm$ of length $m_{a,b}$ in $a$ and $b$, starting with $a$. If $m_{a,b} = m_{b,a} = \infty$, then the relation $\spkl{ab}^{m_{a,b}} = \spkl{ba}^{m_{b,a}}$ means that there is no relation involving $a$ and $b$. The Artin monoid $A_M^+$ of $M$ is the monoid given by the same presentation
$$
  A_M^+ \defeq \spkl{\cS \ \vert \ \spkl{ab}^{m_{a,b}} = \spkl{ba}^{m_{b,a}} \ \forall \ a, b \in \cS}^+.
$$
The semigroup C*-algebras attached to Artin monoids have been studied in \cite{CL02,CL07,LOS}. It is easy to see that $(A_M^+)^* = \gekl{1}$. Moreover, it is shown in \cite{BS72} that $A_M^+$ is right LCM. The main result in \cite{Par} says that $A_M^+$ embeds into $A_M$ via the canonical map. However, it is not clear whether $A_M^+ \subseteq A_M$ satisfies the Toeplitz condition in the sense of \cite{Li13}, or in other words, whether $A_M^+ \subseteq A_M$ is quasi-lattice ordered in the sense of \cite{Nic}. This is only known in special cases, for instance for spherical Artin monoids and groups (also called Artin monoids and groups of finite type) or right-angled Artin monoids and groups, but not in general. Hence we cannot apply the K-theory formula from \cite{CEL2} in general. Nevertheless, Corollary~\ref{cor:LCM} applies and yields the following (with $\cA$ and $A$ as in Corollary~\ref{cor:LCM}):
\setlength{\parindent}{0cm} \setlength{\parskip}{0cm}

\begin{enumerate}
\item[(I)] If $A_M$ satisfies the Baum-Connes conjecture for $\cA$ and $A$, then $K_0(C^*_{\lambda}(A_M^+)) = \Zz [1]_0$ and $K_1(C^*_{\lambda}(A_M^+)) \cong \gekl{0}$.
\item[(II)] If $A_M$ satisfies the strong Baum-Connes conjecture in the sense of \cite[Definition~3.4.17]{CELY}, then the unital embedding $\Cz \into C^*_{\lambda}(A_M^+)$ induces a KK-equivalence between $\Cz$ and $C^*_{\lambda}(A_M^+)$.
\end{enumerate}
\setlength{\parindent}{0cm} \setlength{\parskip}{0.5cm}

\subsection{Semigroup C*-algebras of Baumslag-Solitar monoids}
\label{ss:BS}

Our second example concerns Baumslag-Solitar monoids. Let $k, l \in \Zz$ be non-zero integers. The Baumslag-Solitar group $BS(k,l)$ (see \cite{BS62}) is given by the presentation
$$
  BS(k,l) \defeq \spkl{a,b \ \vert \ ab^k = b^la}.
$$
Baumslag-Solitar monoids are defined in an analogous fashion, but we have to adjust the defining relation in order to avoid inverses. We define the following monoids by presentations:
$$
  BS(k,l)^+ \defeq 
  \bfa  
  \spkl{a,b \ \vert \ ab^k = b^la}^+ & \rm{if} \ k, l > 0.\\
  \spkl{a,b \ \vert \ a = b^lab^{-k}}^+ & \rm{if} \ k<0, \, l > 0.\\
  \spkl{a,b \ \vert \ b^{-l}ab^k = a}^+ & \rm{if} \ k>0, \, l<0.\\
  \spkl{a,b \ \vert \ b^{-l}a = ab^{-k}}^+ & \rm{if} \ k, l< 0.  
  \efa
$$
The semigroup C*-algebras of $BS(k,l)^+$ have been studied in \cite{Spi}. It is easy to see that $(BS(k,l)^+)^* = \gekl{1}$. Using normal forms (see for instance \cite{SW}), it follows that $BS(k,l)^+$ embeds into $BS(k,l)$ via the canonical map. Since $BS(k,l)$ has the Haagerup property by \cite{GJ}, \cite{HK} yields that $BS(k,l)$ satisfies the strong Baum-Connes conjecture. Moreover, it is shown in \cite{Spi} that $BS(k,l)^+$ is right LCM. 

However, for $k<-1, l>0$ (or $k>1, l<0$), $BS(k,l)^+$ does not embed into a group $G$ such that $BS(k,l)^+ \subseteq G$ satisfies the Toeplitz condition. To show this, it suffices to show that $BS(k,l)^+ \subseteq BS(k,l)$ does not satisfy the Toeplitz condition (see \cite[Corollary~5.8.9]{CELY}). The latter claim follows essentially from computations in \cite{Spi}. Indeed, suppose that $BS(k,l)^+ \subseteq BS(k,l)$ satisfies the Toeplitz condition. Let us write $P \defeq BS(k,l)^+$ and $G \defeq BS(k,l)$. Consider the element $g = aba^{-1} \in G$. Since $P$ is right LCM, and because we assume that $P \subseteq G$ satisfies the Toeplitz condition, \cite[Lemma~4.2]{Li13} implies that $gP \cap P = pP$ for some $p \in P$. Let $\#_a$ count the number of $a$s in elements of $P$. If $\#_a(p) = 0$, then $p = b^m$ for some $m \geq 0$. So $b^m = p \in gP$ implies that $b^m = a b a^{-1} x$ for some $x \in P$. Looking at normal forms, we conclude that $x = ay$ for some $y \in P$. But then $b^m = a b y$. Comparing $\#_a$ leads to a contradiction. Now it is easy to see that $ab^n \in gP \cap P$ for all $n \in \Zz$ as $ab^{n-1} \in P$. In particular, $a \in pP$, which implies that $\#_a(p) \leq 1$ and thus $\#_a(p) = 1$. So $p = b^i a b^j$ for some $j \in \Zz$. We can always arrange $0 \leq i < l$. Hence, for every $n \in \Zz$, there exists $x_n \in P$ with $p x_n = a b^n$. Comparing $\#_a$, we obtain that $x_n = b^{k_n}$ for some $k_n \geq 0$. So we have for all $n \in \Zz$ that $b^i a b^{j+k_n} = a b^n$, which implies $i=0$ and $j \leq j+k_n = n$ for all $n \in \Zz$. But this is a contradiction. 
\setlength{\parindent}{0.5cm} \setlength{\parskip}{0cm}

Even worse, it turns out that  $BS(k,l)^+ \subseteq BS(k,l)$ does not even satisfy the weak Toeplitz condition from \cite[Definition~4.5]{CEL2}.
\setlength{\parindent}{0cm} \setlength{\parskip}{0.5cm}

This means that we cannot apply the K-theory formula from \cite{CEL2}. Nevertheless, Corollary~\ref{cor:LCM} allows us to compute K-theory, and we obtain that, for all $k, l \in \Zz \setminus \gekl{0}$, the unital embedding $\Cz \into C^*_{\lambda}(BS(k,l)^+)$ induces a KK-equivalence between $\Cz$ and $C^*_{\lambda}(BS(k,l)^+)$.

\subsection{Semigroup C*-algebras of one-relator monoids}
\label{ss:OneRel}

Our third example is about more classes of one-relator monoids, i.e., monoids of the form $P = \spkl{\cS \ \vert \ u = v}^+$, where $\cS$ is a countable set and $u$, $v$ are finite words in $\cS$. In the following, $\equiv$ stands for equality as finite words, whereas $=$ stands for equality as elements of $P$. We make the following assumptions: First, we always assume that $u \not\equiv \ve \not\equiv v$, where $\ve$ is the empty word. This will imply that $P^* = \gekl{1}$. Secondly, we assume that no $a \in \cS$ is redundant, i.e., for every $a \in \cS$ and $w \in (\cS \setminus \gekl{a})^*$ that $a \neq w$. Thirdly, we always assume that $u \not\equiv v$, and even more, that the first letter of $u$ does not coincide with the first letter of $v$. In this case, if we define the corresponding one-relator group by $G = \spkl{\cS \ \vert \ u = v}$, then $P$ embeds into $G$ via the canonical map. The semigroup C*-algebras of such one-relator monoids have been studied in \cite{LOS}. As observed in \cite[\S~2.1.4]{LOS}, $P$ is right LCM if $\ell^*(u) = \ell^*(v)$ or if $\ell^*(u) < \ell^*(v)$ and there exists $a \in \cS$ with $\ell_a(u) > \ell_a(v)$. Here $\ell^*$ stands for word-length and $\ell_a$ counts how many times $a$ appears. It has been shown in \cite{BBV,Tu,Oyo} that the one-relator group $G$ satisfies the strong Baum-Connes conjecture. Again, although it is not clear whether $P \subseteq G$ satisfies the Toeplitz condition (or equivalently in this case, whether $P \subseteq G$ is quasi-lattice ordered), Corollary~\ref{cor:LCM} nevertheless applies and yields that if $P$ is right LCM, then the unital embedding $\Cz \into C^*_{\lambda}(P)$ induces a KK-equivalence between $\Cz$ and $C^*_{\lambda}(P)$.

Now assume that $\abs{\cS} \geq 3$. Then by \cite[Corollar~3.5]{LOS}, the boundary quotient $\partial C^*_{\lambda}(P)$ (see for instance \cite[\S~5.7]{CELY} for an introduction) is purely infinite simple, and there is an exact sequence $0 \to \cK(\ell^2P) \to C^*_{\lambda}(P) \to \partial C^*_{\lambda}(P) \to 0$ if $3 \leq \abs{\cS} < \infty$ whereas $C^*_{\lambda}(P) = \partial C^*_{\lambda}(P)$ if $\abs{\cS} = \infty$. As explained in \cite[Remark~3.6]{LOS}, our K-theory computation for $C^*_{\lambda}(P)$ then yields the following K-theory formula for $\partial C^*_{\lambda}(P)$:
$$
  (K_0(\partial C^*_{\lambda}(P)), [1]_0, K_1(\partial C^*_{\lambda}(P))) \cong  
  \bfa  
  (\Zz / (\abs{\cS} - 2) \Zz, 1, \gekl{0}) & \rm{if} \ 3 \leq \abs{\cS} < \infty;\\
  (\Zz, 1, \gekl{0}) & \rm{if} \ \abs{\cS} = \infty.  
  \efa
$$
Moreover, \cite[Corollary~3.5]{LOS} implies that $\partial C^*_{\lambda}(P)$ is nuclear if and only if $C^*_{\lambda}(P)$ is nuclear. If that is the case, then our K-theory computations, together with \cite[Chapter~8]{Ror}, imply that
$$
  \partial C^*_{\lambda}(P) \cong  
  \bfa  
  \cO_{\abs{\cS} - 1} & \rm{if} \ 3 \leq \abs{\cS} < \infty;\\
  \cO_{\infty} & \rm{if} \ \abs{\cS} = \infty.  
  \efa
$$
Here $\cO_{\abs{\cS} - 1}$ and $\cO_{\infty}$ denote the corresponding Cuntz algebras. 

In addition, using the notation $E_n^k$ to denote extension algebras of the form $0 \to \cK \to E_n^k \to \cO_n \to 0$ as in \cite[Definition~3.2]{ELR}, \cite[Theorem~3.1]{ELR} implies the following result:
\bcor
\label{cor:OneRelClass}
If $C^*_{\lambda}(P)$ is nuclear, then 
$$
  C^*_{\lambda}(P) \cong  
  \bfa  
  E^{-1}_{\abs{\cS} - 1} & \rm{if} \ 3 \leq \abs{\cS} < \infty;\\
  \cO_{\infty} & \rm{if} \ \abs{\cS} = \infty.  
  \efa
$$
In particular, given two one-relator monoids $P_1 = \spkl{\cS_1 \ \vert \ u_1 = v_1}^+$ and $P_2 = \spkl{\cS_2 \ \vert \ u_2 = v_2}^+$ as above such that their semigroup C*-algebras $C^*_{\lambda}(P_1)$ and $C^*_{\lambda}(P_2)$ are nuclear, we have $C^*_{\lambda}(P_1) \cong C^*_{\lambda}(P_2)$ if and only if $\abs{\cS_1} = \abs{\cS_2}$.
\ecor
Sufficient conditions for nuclearity of $C^*_{\lambda}(P)$ are given in \cite[\S~3]{LOS} (i.e., conditions~(1), (2) and (3) by \cite[Theorem~3.9]{LOS}). Concrete examples where Corollary~\ref{cor:OneRelClass} applies are given in \cite[Example~3.11]{LOS} and look as follows: Let $A$ and $B$ be countable sets with $\abs{A} + \abs{B} \geq 3$ and set $\cS \defeq A \amalg B$. Choose $u \in A^*$ arbitrary and $v \in B^*$ with $\menge{x \in \cS^*}{v \equiv xy \equiv wx \ {\rm for} \ {\rm some} \ \ve \neq w, y \in \cS^*} = \gekl{\ve}$. Then Corollary~\ref{cor:OneRelClass} applies to monoids of the form $\spkl{A \amalg B \, \vert \, u = v}^+$.

\subsection{C*-algebras generated by right regular representations of $ax+b$-type semigroups attached to congruence monoids}
\label{ss:Crho}

Our fourth example class is given by C*-algebras generated by right regular representations of $ax+b$-type semigroups attached to congruence monoids. Let us first introduce the setting. Let $K$ be a number field with ring of algebraic integers $R$. Each fractional ideal $\mfa$ of $K$ can be uniquely written as $\mfa = \prod_{\mfp \in\cP_K} \mfp^{v_{\mfp}(\mfa)}$, where the product runs over the set $\cP_K$ of non-zero prime ideals of $R$ and $v_{\mfp}(\mfa) \in \Zz$ is zero for all but finitely many $\mfp$. Let $\fm = \fm_{\infty} \fm_0$ be a pair consisting of a non-zero ideal $\fm_0$ of $R$ and a collection $\fm_{\infty}$ of real places of $K$. For a real place $w$ of $K$, we write $w \mid \fm_{\infty}$ for $w \in \fm_{\infty}$. Set $R_{\fm} \defeq \menge{a \in R^\times}{v_{\mfp}(a)=0 \text{ for all } \mfp \mid \fm_0}$. Moreover, define $(R/\fm)^* \defeq \rukl{\prod_{w \mid \fm_\infty} \spkl{\pm 1}} \times (R/\fm_0)^*$, and for $a \in R_{\fm}$, let $[a]_{\fm} \defeq (({\rm sign}(w(a)))_{w \mid \fm_\infty},a + \fm_0)\in(R/\fm)^*$. Then $R_{\fm} \to (R/\fm)^*$, $a\mapsto [a]_{\fm}$, is a semigroup homomorphism. If $\Gamma$ is a subgroup of $(R/\fm)^*$, then $R_{\fm,\Gamma} \defeq \menge{a\in R_{\fm}}{[a]_{\fm} \in \Gamma}$ is called a congruence monoid. Let $R\rtimes R_{\fm,\Gamma}$ be the semi-direct product with respect to the action of $R_{\fm,\Gamma}$ on $R$ by multiplication. The semigroup C*-algebra of $R\rtimes R_{\fm,\Gamma}$, more precisely the C*-algebra $C^*_{\lambda}(R\rtimes R_{\fm,\Gamma})$ generated by the left regular representation of $R\rtimes R_{\fm,\Gamma}$, has been studied in \cite{Bru} and \cite{BL}. In particular, in \cite[\S~4.1]{BL}, the K-theory of $C^*_{\lambda}(R\rtimes R_{\fm,\Gamma})$ has been computed. Let us now compute K-theory for the C*-algebra $C^*_{\rho}(R\rtimes R_{\fm,\Gamma})$ generated by the right regular representation of $R\rtimes R_{\fm,\Gamma}$. Since we can identify $C^*_{\rho}(R\rtimes R_{\fm,\Gamma})$ with the C*-algebra $C^*_{\lambda}((R\rtimes R_{\fm,\Gamma})^{\rm op})$ generated by the left regular representation of the opposite semigroup $(R\rtimes R_{\fm,\Gamma})^{\rm op}$, we can apply Corollary~\ref{cor:P} to $(R\rtimes R_{\fm,\Gamma})^{\rm op}$.

Let $K_{\fm, \Gamma} \defeq \menge{x\in K_{\fm}}{[x]_{\fm} \in \Gamma}$. By \cite[Proposition~3.3]{Bru}, the semigroup $R\rtimes R_{\fm,\Gamma}$ is left Ore with group of left quotients equal to $G(R\rtimes R_{\fm,\Gamma}) = (R_{\fm}^{-1} R) \rtimes K_{\fm,\Gamma}$, where $R_{\fm}^{-1} R = \menge{a/b}{a \in R, b \in R_{\fm}} \subseteq K\reg$ is the localization of $R$ at $R_{\fm}$. Hence $P \defeq (R\rtimes R_{\fm,\Gamma})^{\rm op}$ also embeds into $G \defeq (R_{\fm}^{-1} R) \rtimes K_{\fm,\Gamma}$. The group $G$ is solvable, hence amenable, so that it satisfies the strong Baum-Connes conjecture by \cite{HK}. For a non-zero ideal $\mfa$ of $R$, set $\mfa_{\fm,\Gamma} \defeq \mfa \cap R_{\fm,\Gamma}$. Now it follows from \cite[Proposition~6.1]{CEL2} and \cite[Proposition~3.9]{Bru} that $\cJ_P\reg = \menge{R \times \mfa_{\fm,\Gamma}}{(0) \neq \mfa \triangleleft R}$. Furthermore, \cite[Proposition~3.9]{Bru} implies that $P$ satisfies the independence condition. While it is in general not clear whether $P \subseteq G$ satisfies the Toeplitz condition, we can nevertheless compute K-theory of $C^*_{\rho}(R\rtimes R_{\fm,\Gamma}) \cong C^*_{\lambda}(P)$ by applying Corollary~\ref{cor:P}. Note that in general, $P$ is no longer right LCM. In order to present the K-theory formula, let $C_{\fm}^{\bar{\Gamma}} \defeq \cI_{\fm} / K_{\fm,\Gamma}$ be the quotient of the group $\cI_{\fm}$ of fractional ideals of $K$ coprime to $\fm_0$ under the multiplication action of $K_{\fm,\Gamma}$, and given a fractional ideal $\mfa$ of $K$, set $(R:\mfa) \defeq \menge{x \in K}{x \mfa \subseteq R}$. Moreover, given $\fk = [\mfa_{\fk}] \in C_{\fm}^{\bar{\Gamma}}$, set $i_{\fk}: \: C^*_{\rho}((R:\mfa_{\fk}) \rtimes R_{\fm,\Gamma}^*) \to C^*_{\rho}(R\rtimes R_{\fm,\Gamma}), \, \rho_g \ma \bm{1}_{R \times (\mfa_{\fk})_{\fm,\Gamma}} \rho(g)$, where $g \ma \rho_g$ and $p \ma \rho(p)$ are the right regular (anti-)representations of $(R:\mfa_{\fk}) \rtimes R_{\fm,\Gamma}^*$ and $R\rtimes R_{\fm,\Gamma}$, respectively, and $\bm{1}_{R \times (\mfa_{\fk})_{\fm,\Gamma}}$ is the characteristic function on ${R \times (\mfa_{\fk})_{\fm,\Gamma}}$ (which lies in $C^*_{\rho}(R\rtimes R_{\fm,\Gamma})$). We now obtain the following K-theory formula:
\bcor
\label{cor:Kax+b}
$\sum_{\fk \in C_{\fm}^{\bar{\Gamma}}} i_{\fk}$ induces a KK-equivalence between $\bigoplus_{\fk \in C_{\fm}^{\bar{\Gamma}}} C^*_{\rho}((R:\mfa_{\fk}) \rtimes R_{\fm,\Gamma}^*)$ and $C^*_{\rho}(R\rtimes R_{\fm,\Gamma})$.

We have $C^*_{\lambda}(R\rtimes R_{\fm,\Gamma}) \sim_{KK} C^*_{\rho}(R\rtimes R_{\fm,\Gamma})$.
\ecor
\setlength{\parindent}{0cm} \setlength{\parskip}{0cm}

\bproof
For the first claim, it suffices to show that the stabilizer group $G_{R \times \mfa_{\fm,\Gamma}} = \menge{g \in G}{(R \times \mfa_{\fm,\Gamma}) \cdot g = R \rtimes \mfa_{\fm,\Gamma}}$ is given by $G_{R \times \mfa_{\fm,\Gamma}} = (R:\mfa) \rtimes R_{\fm,\Gamma}^*$, for all non-zero ideals $\mfa$ of $R$ coprime to $\fm_0$. \an{$\supseteq$} is clear. To prove \an{$\subseteq$}, take $(b,a) \in G_{R \times \mfa_{\fm,\Gamma}}$. Looking at the multiplicative component, we see that $\mfa_{\fm,\Gamma} \cdot a = \mfa_{\fm,\Gamma}$. As $\mfa_{\fm,\Gamma}$ generates $\mfa$ as an ideal in $R$ by \cite[Lemma~3.8]{Bru}, it follows that $\mfa \cdot a = \mfa$. Thus $a \in R^* \cap K_{\fm,\Gamma} = R_{\fm,\Gamma}^*$. Now looking at the additive component, $(R \rtimes \mfa_{\fm,\Gamma}) \cdot (b,a) = R \rtimes \mfa_{\fm,\Gamma}$ implies that $R + \mfa_{\fm,\Gamma} \cdot b = R$, which is equivalent to $\mfa_{\fm,\Gamma} \cdot b \subseteq R$, which in turn is equivalent to $\mfa \cdot b \subseteq R$ because $\mfa_{\fm,\Gamma} $ generates $\mfa$ as an ideal in $R$ by \cite[Lemma~3.8]{Bru}. Hence $b$ must lie in $(R:\mfa)$, as desired.
\setlength{\parindent}{0cm} \setlength{\parskip}{0.5cm}

For the second claim, just observe that, since $C_{\fm}^{\bar{\Gamma}}$ is a group under multiplication of ideals, the map $[\mfa] \ma [(R:\mfa)]$ defines a bijection on $C_{\fm}^{\bar{\Gamma}}$ because it coincides with the map sending a group element of $C_{\fm}^{\bar{\Gamma}}$ to its inverse. Now our second claim follows from the first claim and \cite[Theorem~4.1]{BL}.
\eproof
\setlength{\parindent}{0cm} \setlength{\parskip}{0.5cm}

In \cite[\S~6]{CEL2} and \cite[\S~4]{Li16}, classes of semigroups where found with the property that their left and right semigroup C*-algebras have the same K-theory, or are even KK-equivalent (see also the discussion in \cite[\S~5.11]{CELY}). Corollary~\ref{cor:Kax+b} identifies more examples with this phenomenon.

\subsection{C*-algebras of inverse semigroups from tilings and point-sets}

As a last class of examples, let us discuss tiling inverse semigroups, point-set inverse semigroups and other related constructions. We refer the reader to \cite{Kel97a,Kel97b,KL00,KL04} for more details. For our K-theory computations for the reduced C*-algebras of the inverse semigroups, it turns out that it is very helpful to have flexibility in choosing the target group of the idempotent pure partial homomorphism on our inverse semigroup. There is always a universal group (see \cite{KL04}), but this group can be difficult to determine, so that it is not so easy to check that this group satisfies the Baum-Connes conjecture with the coefficients of interest. However, all we need in order to apply Theorem~\ref{thm:InvSgp} is to find some idempotent pure partial homomorphism whose target group has the desired properties. This is much easier to achieve, as we will see in the examples below.

Let us start with tiling inverse semigroups. A tile is a subset of $\Rz^n$ which is homeomorphic to a closed ball in $\Rz^n$. A partial tiling is a collection of tiles with pairwise disjoint interiors. The support of a partial tiling is the union of its tiles. A tiling is a partial tiling whose support is all of $\Rz^n$. A patch is a finite partial tiling. Let $\cT$ be a tiling. Let $\cP$ be the set of subpatches of $\cT$. Define an equivalence relation on triples of the form $(a,P,b)$ with $P \in \cP$, $a, b \in P$ by setting $(a,P,b) \sim (c,Q,d)$ if and only if there exists $x \in \Rz^n$ such that $a+x=c$, $b+x=d$ and $P+x=Q$. Let $[\cdot]$ denote equivalence classes with respect to $\sim$. Then $\Gamma(\cT) \defeq \menge{[a,P,b]}{P \in \cP, \, a,b \in P} \cup \gekl{0}$ becomes an inverse semigroup under the multiplication $[a,P,b] \cdot [c,Q,d] \defeq [a+x,(P+x) \cup (Q+y),d+y]$ if there exist $x,y \in \Rz^n$ such that $P+x$ and $Q+y$ are subpatches of $\cT$ and $b+x = c+y$; otherwise define $[a,P,b] \cdot [c,Q,d] \defeq 0$. It is easy to see that $[a,P,b]^{-1} = [b,P,a]$. We call $\Gamma(\cT)$ the tiling inverse semigroup of $\cT$. If we replace $\cP$ by the set $\cP_{\rm conn}$ of subpatches of $\cT$ with connected support and perform the above construction, then we obtain the connected tiling inverse semigroup $S(\cT)$. 

For each tile $t \in \cT$, let us choose a point $p(t)$ in the interior of $t$ ($p(t)$ is called the puncture of $t$) such that if for $x \in \Rz^n$, both $t$ and $t+x$ are tiles in $\cT$, then $p(t+x) = p(t)+x$. Now let $G \defeq \spkl{\menge{p(t) - p(t')}{t,t' \in \cT}} \subseteq \Rz^n$ be the additive subgroup of $\Rz^n$ generated by $p(t) - p(t')$ for $t, t' \in \cT$. Since $\cT$ is countable, $G$ is a countable group. It is straightforward to check that $\sigma: \: \Gamma(\cT)\reg \to G, \, [a,P,b] \ma p(a) - p(b)$ defines an idempotent pure partial homomorphism. Similarly, the restriction of $\sigma$ to $S(\cT)\reg$ defines an idempotent pure partial homomorphism on $S(\cT)\reg$. The group $G$ is abelian, hence satisfies the strong Baum-Connes conjecture by \cite{HK}. Thus we can apply Theorem~\ref{thm:InvSgp} to compute K-theory for $C^*_{\lambda}(\Gamma(\cT))$ and $C^*_{\lambda}(S(\cT))$. To present the K-theory formula, we introduce the equivalence relation $\approx$ on $\cP$ and $\cP_{\rm conn}$ by setting $P \approx Q$ if and only if there exists $x \in \Rz^n$ with $Q = P+x$. Given $P \in \cP$, choose $a \in P$ and denote by $i_P$ the homomorphism $\Cz \to C^*_{\lambda}(\Gamma(\cT))$ (or $\Cz \to C^*_{\lambda}(S(\cT))$) sending $1 \in \Cz$ to $[a,P,a]$. \cite[Lemmas~6.1 and 6.2]{Nor15} together with Theorem~\ref{thm:InvSgp} now yield the following:
\bcor
$\sum_{[P] \in \cP/\approx} i_P$ induces a KK-equivalence between $\bigoplus_{[P] \in \cP/\approx} \Cz$ and $C^*_{\lambda}(\Gamma(\cT))$. 

$\sum_{[P] \in \cP_{\rm conn}/\approx} i_P$ induces a KK-equivalence between $\bigoplus_{[P] \in \cP_{\rm conn}/\approx} \Cz$ and $C^*_{\lambda}(S(\cT))$. 
\ecor 
\setlength{\parindent}{0cm} \setlength{\parskip}{0cm}

This generalizes \cite[Proposition~6.3]{Nor15}. The reason we can now cover all tiling inverse semigroups is that we no longer need the (much) stronger condition in \cite{Nor15} that our inverse semigroups have to be $0$-$F$-inverse semigroups and must admit a partial homomorphism to a group which is injective on maximal elements.
\setlength{\parindent}{0cm} \setlength{\parskip}{0.5cm}

Let us now discuss point-set inverse semigroups. We start with a countable subset $\cD \subseteq \Rz^n$. Let $\cP$ be the set of finite subsets of $\cD$. Define an equivalence relation on triples of the form $(a,P,b)$ with $P \in \cP$, $a, b \in P$ by setting $(a,P,b) \sim (c,Q,d)$ if and only if there exists $x \in \Rz^n$ such that $a+x=c$, $b+x=d$ and $P+x=Q$. Let $[\cdot]$ denote equivalence classes with respect to $\sim$. Then $\Gamma(\cD) \defeq \menge{[a,P,b]}{P \in \cP, \, a,b \in P} \cup \gekl{0}$ becomes an inverse semigroup under the multiplication $[a,P,b] \cdot [c,Q,d] \defeq [a+x,(P+x) \cup (Q+y),d+y]$ if there exist $x,y \in \Rz^n$ such that $P+x$ and $Q+y$ are finite subsets of $\cD$ and $b+x = c+y$; otherwise define $[a,P,b] \cdot [c,Q,d] \defeq 0$. We call $\Gamma(\cD)$ the point-set inverse semigroup of $\cD$.

Now let $G \defeq \spkl{\menge{d - d'}{d,d' \in \cD}} \subseteq \Rz^n$ be the additive subgroup of $\Rz^n$ generated by differences of elements of $\cD$. Since $\cD$ is countable, $G$ is a countable group. Moreover, $\sigma: \: \Gamma(\cD)\reg \to G, \, [a,P,b] \ma a - b$ defines an idempotent pure partial homomorphism. As $G$ is abelian, it satisfies the strong Baum-Connes conjecture by \cite{HK}. Thus we can apply Theorem~\ref{thm:InvSgp} to compute K-theory for $C^*_{\lambda}(\Gamma(\cD))$. As above, we introduce the equivalence relation $\approx$ on $\cP$ by setting $P \approx Q$ if and only if there exists $x \in \Rz^n$ with $Q = P+x$. Given $P \in \cP$, choose $a \in P$ and denote by $i_P$ the homomorphism $\Cz \to C^*_{\lambda}(\Gamma(\cD))$ sending $1 \in \Cz$ to $[a,P,a]$. The analogues of \cite[Lemmas~6.1 and 6.2]{Nor15} together with Theorem~\ref{thm:InvSgp} now yield the following:
\bcor
$\sum_{[P] \in \cP/\approx} i_P$ induces a KK-equivalence between $\bigoplus_{[P] \in \cP/\approx} \Cz$ and $C^*_{\lambda}(\Gamma(\cD))$. 
\ecor

Finally, let us discuss inverse semigroups of the form $\Gamma(X,G,H)$ from \cite[Example~2.1.1~(iii)]{KL04}, which are constructed as follows: Let $H$ be a group, $G$ a subgroup of $H$ and $X$ a subset of $H$ with $1 \in X$. Let $\cP$ be the set of finite intersections of subsets of $H$ of the form $\menge{gX}{g \in G, \, 1 \in gX}$. Define an equivalence relation on triples of the form $(a,P,b)$ with $P \in \cP$, $a, b \in P$ by setting $(a,P,b) \sim (c,Q,d)$ if and only if there exists $g \in G$ such that $g \cdot a = c$, $g \cdot b = d$ and $g \cdot P = Q$. Let $[\cdot]$ denote equivalence classes with respect to $\sim$. Then $\Gamma(X,G,H) \defeq \menge{[a,P,b]}{P \in \cP, \, a,b \in P} \cup \gekl{0}$ becomes an inverse semigroup under the multiplication $[a,P,b] \cdot [c,Q,d] \defeq [g \cdot a, (g\cdot P) \cap (h \cdot Q),h \cdot d]$ if there exist $g,h \in G$ such that $g \cdot P, \, h \cdot Q \in \cP$ and $g \cdot b = h \cdot c$; otherwise define $[a,P,b] \cdot [c,Q,d] \defeq 0$. 

It is straightforward to check that $\sigma: \: \Gamma(X,G,H)\reg \to G, \, [a,P,b] \ma ab^{-1}$ defines an idempotent pure partial homomorphism. So we can apply Theorem~\ref{thm:InvSgp} if $G$ is a countable group satisfying the Baum-Connes conjecture with the relevant coefficients or in its strong form. As above, we introduce the equivalence relation $\approx$ on $\cP$ by setting $P \approx Q$ if and only if there exists $g \in G$ with $Q = g \cdot P$. Given $P \in \cP$, choose $a \in P$ and denote by $i_P$ the homomorphism $\Cz \to C^*_{\lambda}(\Gamma(X,G,H))$ sending $1 \in \Cz$ to $[a,P,a]$. Let $\cA$ and $A$ be the C*-algebras constructed in \S~\ref{ss:MoritaEnv} and \S~\ref{ss:DiscreteV} for the inverse semigroup $\Gamma(X,G,H)$. The analogues of \cite[Lemmas~6.1 and 6.2]{Nor15} together with Theorem~\ref{thm:InvSgp} now yield the following:
\bcor
$ $
\setlength{\parindent}{0cm} \setlength{\parskip}{0cm}

\begin{enumerate}
\item[(I)] If $G$ is countable and satisfies the Baum-Connes conjecture for $\cA$ and $A$, then $\sum_{[P] \in \cP/\approx} i_P$ induces a K-theory isomorphism $\bigoplus_{[P] \in \cP/\approx} K_*(\Cz) \cong K_*(C^*_{\lambda}(\Gamma(X,G,H)))$.
\item[(II)] If $G$ is countable and satisfies the strong Baum-Connes conjecture in the sense of \cite[Definition~3.4.17]{CELY}, then $\sum_{[P] \in \cP/\approx} i_P$ induces a KK-equivalence between $\bigoplus_{[P] \in \cP/\approx} \Cz$ and $C^*_{\lambda}(\Gamma(X,G,H))$.
\end{enumerate}
\ecor
\setlength{\parindent}{0cm} \setlength{\parskip}{0.5cm}

\end{document}